\documentclass[11pt]{article}
\usepackage{hyperref}
 \usepackage{benstyle}
 \usepackage{cite}
 
\usepackage{enumitem}
\usepackage{authblk}

\newcommand*{\alphanew}{\tau} 
\newcommand*{\betanew}{\omega}

\title{Towards Sharp Minimax Risk Bounds for Operator Learning}

\author[1]{Ben Adcock}
\author[2,3]{Gregor Maier}
\author[4]{Rahul Parhi}

\affil[1]{Department of Mathematics, Simon Fraser University, Burnaby~BC, Canada}
\affil[2]{Institute for Numerical Simulation, University of Bonn, Bonn, Germany}
\affil[3]{Fraunhofer Institute for Algorithms and Scientific Computing (SCAI), Sankt Augustin, Germany}
\affil[4]{Department of Electrical and Computer Engineering, University of California, San Diego, La Jolla, CA, USA}

\begin{document}

\maketitle

\begin{abstract}
We develop a minimax theory for operator learning, where the goal is to estimate an unknown operator between separable Hilbert spaces from finitely many noisy input-output samples. For uniformly bounded Lipschitz operators, we prove information-theoretic lower bounds together with matching or near-matching upper bounds, covering both fixed and random designs under Hilbert-valued Gaussian noise and Gaussian white noise errors. The rates are controlled by the spectrum of the covariance operator of the measure that defines the error metric. Our setup is very general and allows for measures with unbounded support. A key implication is a curse of sample complexity, which shows that the minimax risk for generic Lipschitz operators cannot decay at any algebraic rate in the sample size. We obtain sharp characterizations when the covariance spectrum decays exponentially and provide general upper and lower bounds in slower-decay regimes. Finally, we show that assuming higher regularity, i.e., H\"older smoothness, does not improve minimax rates over the Lipschitz case, up to potential constants. Thus, we show that learning operators of any finite regularity necessarily suffers a curse of sample complexity.
\end{abstract}


\pbk
\noindent \textbf{Corresponding author:} \url{ben_adcock@sfu.ca} (Ben Adcock)

\section{Introduction}

A new paradigm in machine learning for scientific computing is focused on designing learning algorithms and methods for continuum problems. This paradigm is referred to as \emph{operator learning} and has received considerable interest in the last few years~\cite{boulle2024mathematical,anandkumar2020neural,herde2024poseidon,hesthaven2018non,kovachki2024operator,lanthaler2024operator,lanthaler2023operator,li2021fourier,lu2021learning,raonic2023convolutional,nelsen2024operator}. The basic task may be posed as learning a map between infinite-dimensional function spaces, i.e., learning an operator
\bes{
    F: \cX \to \cY,
}
where, for example, $\cX$ and $\cY$ are real, separable Hilbert spaces. Operator learning naturally arises in many scientific problems where one wants to learn how a continuum model, often described by partial differential equations (PDEs), maps inputs, such as parameters or boundary conditions, to outputs, such as states or observables. A prototypical example to keep in mind is learning parameter-to-solution maps of parametric PDEs \cite{adcock2022sparse,cohen2015approximation,adcock2024learning}. In contrast to more classical surrogate modeling, which typically focuses on learning finite-dimensional parameter-to-solution maps for some \emph{fixed discretization}, operator learning directly aims to learn/approximate the continuum map $F: \cX \to \cY$ itself. Thus, the inputs and outputs are functions (not vectors) and the goal is to directly design \emph{discretization-invariant} methods~\cite{kovachki2024operator,boulle2024mathematical}.

From a statistical perspective, this naturally leads to a nonparametric regression problem in which both the object of interest (the operator) and the observations (finite number of noisy samples) are infinite-dimensional. Suppose that we observe the noisy input-output pairs
\bes{
    \{(X_i, Y_i)\}_{i=1}^m, \quad Y_i = F(X_i) + \sigma E_i, \quad i = 1, \ldots, m,
}
where the design points $\{X_i\}_{i=1}^m$ are either \emph{fixed} or \emph{random} elements of $\cX$, the noise terms $E_i$ are i.i.d.\ noise terms which may or may not take values in $\cY$ (e.g., our setup can account for Gaussian white noise that almost surely does not take values in $\cY$), 
 $\sigma > 0$ is the noise level, and $F$ belongs in some prescribed model class $\cF$ of operators. The goal is then to use the data $\{(X_i, Y_i)\}_{i=1}^m$ to design an \emph{estimator} $\widehat{F} := \widehat{F}_m$ that is reasonably close to $F$. The performance of an estimator is quantified via the \emph{minimax risk} 
\be{
    \inf_{\widehat{F}} \sup_{F \in \cF} \bbE \left [ \nm{F - \widehat{F}}_{L^p_{\mu}(\cX ; \cY)} \right ], \quad 1 \leq p < \infty. \label{eq:minimax-intro}
}
The precise dependence of \eqref{eq:minimax-intro} on $m$ is referred to as the \emph{minimax rate}.
In finite-dimensional nonparametric regression (i.e., $\cX$ and $\cY$ are finite-dimensional Euclidean spaces), the optimal behavior of \eqref{eq:minimax-intro} is well-understood. Indeed, for many finite-dimensional model classes, e.g., the unit ball of the Besov space $B^s_{p,q}[-1,1]^d$ that satisfies the condition $s > d/p$ (so that $B^s_{p,q}[-1,1]^d$ compactly embeds into $C[-1,1]^d$), the minimax rate has been exactly determined~\cite{delyon1996minimax,devore2025optimal,donoho1994ideal,donoho1998minimax,tsybakov2008introduction,vandeGeerBook,WainwrightBook}. For operator learning, however, the situation is much less clear (see Section~\ref{sec:contrib-related} for what is currently known). The infinite-dimensional nature of the input and output spaces is what makes this problem to be intrinsically challenging. This motivates a very fundamental statistical question that has been largely unanswered:

\begin{center}
  \emph{Given a model class $\cF$ of operators $F : \cX \to \cY$ and the error metric of the $L^p_\mu(\cX; \cY)$-norm,
  what is the optimal rate at which the risk can decay as a function of $m$?}
\end{center}
In this paper, we make a step towards closing this gap for one of the most natural and widely studied operator classes, the class of \emph{uniformly bounded Lipschitz operators}. We focus on the nonparametric regression setting where the noise is either Hilbert-valued Gaussian noise (trace-class covariance operator) or Gaussian white noise. For this model class, we study the minimax risk for both fixed and random designs and quantify the decay rates as a function of the eigenvalues $\{\lambda_i\}_{i \geq 1}$ of the covariance operator of the measure $\mu$. Our main results give general upper and lower bounds on this quantity for classes of bounded Lipschitz operators, and we derive sharp rates in several important regimes. In particular, we remark that our setting is very general as $\mu$ could have bounded or unbounded support and includes practical cases such as the uniform measure on a compact set or an (unbounded support) Gaussian measure. This generality is in stark contrast with most works in minimax estimation. Finally, we show that imposing higher smoothness does not help: for operators of higher-order H\"older smoothness, the minimax rates do not improve over those of the Lipschitz case, except for potential constants.

At a qualitative level, our findings show that operator learning is subject to a \emph{curse of sample complexity}: Regardless of how fast the eigenvalues decay, the minimax risk cannot decrease at an algebraic rate in $m$. 
In particular, in the case of algebraically-decaying eigenvalues, $\lambda_i = i^{-\alphanew}$, we establish general upper and lower bounds for the minimax risk, with the former being algebraically-decaying in $\log(m)$ and the latter exponentially-decaying in $\sqrt{\log(m)}$ (note that this subalgebraic in $m$, as must be the case). For exponentially decaying eigenvalues, $\lambda_i = \exp(-\alphanew i^\betanew)$, we obtain essentially matching upper and lower bounds for the minimax risk, leading to a precise characterization of the minimax rate in terms
of $m$. Overall, these results make precise, in a minimax sense, the curse of sample complexity that has been identified in previous data-complexity and approximation-theoretic analyses of operator learning, which we discuss below.

\subsection{Motivations}

There are multiple motivations for developing a minimax theory of operator learning for Lipschitz classes. For example, Lipschitz operators provide a flexible and natural model class for many operator learning problems arising in applications. Maps defined by solution operators of PDEs, parameter-to-state maps, and parameter-to-observable maps often enjoy a Lipschitz (or locally Lipschitz) dependence on the input in appropriate norms. Such Lipschitz bounds play a central role, for example, in stability analyses of PDEs and in the well-posedness theory for inverse problems. Furthermore, it is well-known that solution operators of elliptic variational inequalities, which arise, for instance, in obstacle and optimal control problems, are Lipschitz continuous but generally not Gâteaux differentiable~\cite{christof2018sensitivity,mignot1976controle,christof2020differential}. It is therefore natural to ask how difficult it is, from a statistical point of view, to learn a generic Lipschitz operator from noisy input-output data. 

There has also been recent work in operator learning on revealing various notions of ``complexity'' for operator classes, including \emph{parametric complexity} (the number of parameters required by a neural operator architecture to achieve a given approximation error)~\cite{lanthaler2025parametric,schwab2025deep} 
and \emph{data complexity} (the number of samples required by a particular estimation procedure to achieve a prescribed accuracy)~\cite{adcock2024optimalb,adcock2025sample,kovachki2024data}. Approximation-theoretic results for neural operators have shown that approximating general Lipschitz operators may require an exponentially large number of parameters in the target accuracy. Data-complexity analyses, on the other hand, suggest that sample requirements for operator learning can be extremely unfavorable for Lipschitz operators~\cite{adcock2025sample,kovachki2024data,subedi2025operator}. A minimax perspective complements these analyses by asking a more fundamental question: Even if we ignore computational and architecture choices, what is the best sample complexity that \emph{any} estimator can achieve over a given operator class from a set of $m$ noisy observations?

Thus, the results of this paper contribute to this agenda by providing the first, to the best of our knowledge, general minimax bounds (upper and lower) for bounded Lipschitz operators and by identifying regimes where these bounds are sharp. Along the way, our analysis clarifies how the difficulty of operator learning depends on the geometry induced by the underlying measure $\mu$ (encoded in the eigenvalues $\{\lambda_i\}_{i\geq 1}$). 
Our analysis also reveals many nuances when studying statistical estimation where the inputs and outputs are both infinite dimensional when compared to the classical statistical paradigm of finite-dimensional inputs and outputs, which we recover as a special case (see Section~\ref{sec:finite-dimensional-rates}). 

\subsection{Implications and related work} \label{sec:contrib-related}

This paper makes several contributions that are statistical in nature. Although operator learning has been studied extensively from approximation-theoretic and algorithmic perspectives, there is comparatively little work that characterizes its fundamental statistical difficulty. We develop a minimax theory for learning uniformly bounded Lipschitz operators from noisy input-output samples, providing both information-theoretic lower bounds and matching (or near-matching) upper bounds across fixed and random designs and across two canonical infinite-dimensional noise models. 
The arguments draw together ideas that are typically disparate: Operator learning, minimax estimation under general (possibly unbounded) measures, nonparametric functional data analysis, and estimation of objects with Hilbert-valued outputs. Thus, the results of this paper can also be viewed as yielding a single framework that subsumes several previously studied settings as special cases. We outline these previous settings in the context of our results in the sequel.

\paragraph{Data complexity and minimax rates for operator learning.}
There has recently been a large interest in studying various statistical and data-complexity properties for operator learning problems~\cite{de2023convergence,jin2023minimax,kovachki2024data,mhaskar1997neural,subedi2024online,subedi2025operator}. Outside of a few notable examples~\cite{liu2024deep,reinhardt2024statistical}, there is almost no work on minimax estimation for operator learning.  The authors of~\cite{liu2024deep} analyze ERM solutions over deep network classes for estimating Lipschitz operators between separable Hilbert spaces and derive nonasymptotic upper bounds for the risk. They note that sharp determination of minimax rates (and, in particular, matching lower bounds) is an open problem. On the other hand, the authors of~\cite{reinhardt2024statistical} develop a general statistical theory for determining upper bounds on operators using tools from empirical-process theory. As a special case, they obtain algebraic rates for learning holomorphic operators, but note that lower bounds are unknown in the infinite-dimensional setting. In contrast, we study the minimax risk over uniformly bounded Lipschitz operators and smoother counterparts and prove information-theoretic lower bounds together with matching (or near-matching) upper bounds, implying in particular that the minimax risk cannot decay algebraically in the sample size for generic opertors. We find sharp rates in regimes such as exponential eigenvalue decay.  

\paragraph{Nonparametric estimation beyond bounded domains.}
While classical minimax theory is often stated for compact domains, our framework is formulated for error measured with respect to a general measure $\mu$ (possibly unbounded support) and recovers the usual bounded-domain Lipschitz minimax rate as a special case (e.g., $\mu$ uniform on $[-1,1]^d$ yields $m^{-\frac{1}{2+d}}$, cf. Section~\ref{sec:finite-dimensional-rates}). Recently, a finite-dimensional counterpart on unbounded supports appeared in~\cite{bizeul2025entropy}. There, the authors study scalar regression of $1$-Lipschitz functions $f:\bbR^d\to\bbR$ under isotropic log-concave $\mu$ and analyze various estimators. In contrast, we consider the operator learning problem and characterize the minimax rate in terms of the eigenvalues $\{\lambda_i\}_{i\ge 1}$, giving information-theoretic lower bounds and matching upper bounds in infinite-dimensional settings.

\paragraph{Nonparametric functional data analysis.} 
A close analogue of operator learning is nonparametric functional regression, where one estimates a regression operator from functional designs and possibly Hilbert-valued responses; see, e.g., the book~\cite{ferraty2006nonparametric} and references therein. In this literature, rates for kernel/Nadaraya--Watson and $k$-nearest-neighbor methods are typically controlled by local small-ball probabilities. For instance, the authors of~\cite{burba2009k,lian2011convergence} prove almost-sure convergence rates (only upper bounds) for Nadaraya--Watson and nearest-neighbor estimators with scalar-valued~\cite{burba2009k} and Hilbert-valued~\cite{lian2011convergence} observations (i.e., the noise is $\cY$-valued, meaning that Gaussian white noise is excluded) under boundedness and H\"older/Lipschitz regularity assumptions on the regression operator. However, sharpness of their bounds, i.e., characterization of the minimax rate, was not addressed. In the opposite direction, the authors of~\cite{mas2012lower} derive minimax lower bounds for estimating functionals $F: \cX \to \bbR$. They show that, in general, for infinite-dimensional settings, the minimax risk at \emph{a single point} is necessarily subalgebraic in the number of data points $m$, but, as above, the sharpness of their bounds was not addressed. Our results, on the other hand, provide \emph{Hilbert-valued} and \emph{global-risk} counterparts of these results with matching (or near-matching) upper and lower bounds in several important regimes. 

\paragraph{Estimation of objects with Hilbert-valued outputs.}
A recurring technical issue in infinite-dimensional-output regression/operator learning is the noise model. Many analyses assume that the noise is $\cY$-valued (e.g., sub-Gaussian or trace-class Gaussian)~\cite{lian2011convergence}. This ensures that the observations satisfy $Y_i \in \cY$ almost surely, in which case least-squares objectives are well-defined. Unfortunately, this assumption does not include the Gaussian white noise errors as soon as $\dim(\cY) = \infty$. As a result, naive least-squares objectives are not well-defined, and so one must work with alternative formulations. In this paper we treat both of these canonical regimes (Assumptions~\ref{ass:hilbert-noise} and~\ref{ass:white-noise}). In the case of Gaussian white noise errors we take inspiration from~\cite{de2023convergence,herrmann2024neural,nickl2020convergence,reinhardt2024statistical} and consider Hilbert scales to ``tame'' the white noise, cf.~\eqref{Yt-def}. Correspondingly, our results allow for a unified minimax theory for Hilbert-valued errors and white noise errors.

\subsection{Road map and paper organization}
Section~\ref{sec:main} introduces our main results in the Lipschitz case, and, in particular, states our main bounds across a range of regimes for the eigenvalue decay. The key technical aspects are the general information-theoretic lower bound proved in Section~\ref{sec:lower} and the complementary upper bound proved in Section~\ref{sec:upper}. Section~\ref{s:proofs-of-main-res} includes the proofs of these main results. In Section~\ref{sec:further}, we present further discussion, including the recovery of classical finite-dimensional minimax rates, a sharp analysis in the double-exponential eigenvalue-decay regime, and show that higher-order (H\"older) regularity does not improve the minimax rates. We conclude in Section~\ref{sec:conclusion} with a number of open problems.

\section{Main results and discussion} \label{sec:main}

In this section, we setup the problem, including key assumptions, then present our main results.

\subsection{Setup and assumptions}

We consider operators $F : \cX \rightarrow \cY$, where $\cX$, $\cY$ are real, separable Hilbert spaces. We are generally interested in the case where both $\cX$ and $\cY$ are infinite-dimensional---however, our analysis also covers the finite-dimensional setting.
We equip $\cX$ with a probability measure $\mu$ satisfying the following assumptions:

\begin{assumption}[Assumptions on $\mu$]
\label{ass:mu}
The probability measure $\mu$ satisfies:
\begin{enumerate}[label=(\roman*)]
\item[(i)] $\mu$ has mean zero and finite second moments, i.e. $\int_{\cX} \nm{X}^2_{\cX} \D \mu(X) < \infty$.
\item[(ii)] Let $\lambda_1 \geq \lambda_2 \geq \cdots > 0$ 
and $\{ \phi_i \}_{i} \subset \cX$ be the eigenvalues and orthonormal basis of eigenvectors of the covariance operator of $\mu$ (which exist due to (i)). Then the real-valued random variables defined by $\xi_i = \ip{X}{\phi_i}_{\cX} / \sqrt{\lambda_i}$ for $X \sim \mu$ are independent.
\item[(iii-L)] Each $\xi_i$ is absolutely continuous with respect to the Lebesgue measure $\lambda$ and there exist constants $a,b > 0$ such that the densities $\nu_i = \D \xi_i / \D \lambda$ satisfy
\bes{
\nu_i(x) \geq b,\quad \forall |x| \leq a,\ i = 1,2,\ldots .
}
Note that $2 a b \leq 1$ and this bound is tight. 
We will, for convenience, define $\iota = 2 a b \leq 1$.
\item[(iii-U)] Suppose the error is measured in the $L^p$-norm (see below). Then the $\xi_i$ have uniformly bounded $p'$th moments, where $p' = \max \{ p ,2 \}$, i.e., $\sup_{i} [ \bbE | \xi_i |^{p'} ]^{\frac{1}{p'}} < \infty$.
\end{enumerate}
\end{assumption}
Here (iii-L) denotes the assumption required for the lower bound on the minimax error and (iii-U) denotes the assumption required for the upper bound. 

Throughout, we consider noisy evaluations of $F$ of the form
\be{
\label{F-evals}
Y_i = F(X_i) + \sigma E_i,\quad i = 1,\ldots,m,
}
where the $E_i$ are drawn independently from a suitable noise distribution $\varepsilon$ (see below). The design points $\cD = \{X_1,\ldots,X_m \} \subset \cX$ 
are either fixed, but arbitrary or chosen randomly. In the latter case, where assume that $X_1,\ldots,X_m \sim_{\mathrm{i.i.d.}} \nu$, where $\nu$ is a probability measure on $\cX$, and the $X_i$ and $E_i$ are independent. We consider the following two noise models.

\begin{assumption}[Hilbert-valued Gaussian noise]\label{ass:hilbert-noise}
$E_1,\ldots,E_m \sim_{\mathrm{i.i.d.}} \varepsilon$, where $\varepsilon$ is a mean-zero Gaussian measure on $\cY$ with covariance operator $\Upsilon : \cY \rightarrow \cY$ with $\mathrm{Tr}(\Upsilon) = 1$ and eigenvalues $\upsilon_1 \geq \upsilon_2 \geq \cdots \geq 0$.
\end{assumption}
Note that the normalization $\mathrm{Tr}(\Upsilon) = 1$ is simply for convenience, and can obviously be removed by adjusting $\sigma$ in \ef{F-evals}.

\begin{assumption}[White noise]\label{ass:white-noise}
$E_1,\ldots,E_m \sim_{\mathrm{i.i.d.}} \varepsilon$, where $\varepsilon$ is Gaussian white noise on $\cY$.
\end{assumption}
Let $1 \leq p \leq \infty$ and $\cF$ be a class of operators $F : \cX \rightarrow \cY$. In the fixed design case, we define the minimax risk
\be{
\label{minimax-fixed}
\cM_{m}(\cF  ; L^p_{\mu}) =\cM^{\mathsf{fix}}_{m}(\cF  ; L^p_{\mu}) = \inf_\cD\inf_{\widehat{F}} \sup_{F \in \cF} \bbE \left [ \nm{F - \widehat{F}}_{L^p_{\mu}(\cX ; \cY)} \right ],
}
where the inner infima are taken over all estimators $\widehat{F}$, i.e., mappings from the evaluations \ef{F-evals} to $L^p_{\mu}(\cX ; \cY)$, and designs $\cD = \{X_1,\ldots,X_m\}$, and the expectation is with respect to the $E_i$. 
In the random design case, we define
\be{
\label{minimax-random}
\cM_{m}(\cF  ; L^p_{\mu}) = \cM^{\mathsf{rand}}_{m}(\cF  ; L^p_{\mu}) = \inf_{\nu} \inf_{\widehat{F}} \sup_{F \in \cF} \bbE \left[ \nm{F - \widehat{F}}_{L^p_{\mu}(\cX ; \cY)} \right ],
}
where the first infimum is taken over all probability measures $\nu$ on $\cX$ and the expectation is now taken with respect to both the $X_1,\ldots,X_m \sim_{\mathrm{i.i.d.}} \nu$ and the $E_i$. 
Here and throughout, $L^p_{\mu}(\cX ; \cY)$ is the $L^p$-Bochner space, with norm
\bes{
\nm{F}_{L^p_{\mu}(\cX ; \cY)} = \left ( \int_{\cX} \nm{F(X)}^p_{\cY} \D \mu(X) \right )^{1/p}
}
for $1 \leq p < \infty$ and
\bes{
\nm{F}_{L^{\infty}_{\mu}(\cX ; \cY)} = \esssup_{X \sim \mu} \nm{F(X)}_{\cY}
}
for $p = \infty$. 
In our main results, it is convenient to consider the negative log of the minimax risk. Accordingly, we define
\be{
\label{neg-log-minimax}
\cL_m(\cF ; L^p_{\mu}) = - \log(\cM_m(\cF ; L^p_{\mu})).
}
We also require some notation for asymptotic behaviour. Given parameters $m$ and $r$ and functions $A(m,r),B(m,r) > 0$, we write
\be{
\label{lesssim-notation}
A(m,r) \lesssim_r B(m,r),\quad m \rightarrow \infty,
}
if there is a constant $c = c(r) > 0$ depending on $r$ only such that $A(m) \leq c(r) B(m)$ for all sufficiently large $m$. We define $A(m,r) \gtrsim_r B(m,r)$ analogously, and write $A(m,r) \asymp_r B(m,r)$ if $A(m,r) \lesssim_r B(m,r)$ and $B(m,r) \lesssim_r A(m,r)$. We will, on occasion, also write $A(m,r) = \cO_r(B(m,r))$ with the same meaning as \ef{lesssim-notation}.

\subsection{Main results}\label{ss:main-result}

Our main results consider classes of uniformly bounded Lipschitz operators. Let $B , L \geq 0$. Then we define 
\bes{
\cF_{B,L} = \left \{ F: \cX \to \cY : \nm{F(X)}_{\cY} \leq B,\ \nm{F(X) - F(X')}_{\cY} \leq L \nm{X - X'}_{\cX},\ \forall X,X' \in \cX \right \}.
}
In order to establish results in the case of Gaussian white noise, we restrict to operators taking values in certain smoothness spaces, which are subspaces of $\cY$ equipped with a scale parameter $t \geq 0$. Following \cite{herrmann2024neural,reinhardt2024statistical}, let $\{ \psi_i \}_{i} \subset \cY$ be an arbitrary orthonormal basis and $(w_i)_{i \in \bbN}$ be a nondecreasing sequence of positive weights. For $t \geq 0$, we define the subspaces
\be{
\label{Yt-def}
\cY^t = \left \{ Y \in \cY : \nm{Y}^2_{\cY^t} : = \sum_{i} w^{2t}_i | \ip{Y}{\psi_i}_{\cY} |^2 < \infty \right \}.
}
Note that $\cY^t$ forms a Hilbert space with the resulting inner product.
We will primarily consider the case $w_i = i$ in what follows. 
As we discuss in Remark \ref{wr-choice}, other decay types could readily be analyzed.
Note that when $\cY = L^2(D)$ for some domain $D \subseteq \bbR^k$ then the spaces $\cY^t$ often correspond to Sobolev spaces, depending on the choice of basis $\{ \psi_i \}_i$. See \cite{herrmann2024neural} and references therein for further discussion.
With this in hand, we define
\bes{
\cF^t_{B,L} =  \left \{ F: \cX \to \cY : \nm{F(X)}_{\cY^t} \leq B,\ \nm{F(X) - F(X')}_{\cY} \leq L \nm{X - X'}_{\cX},\ \forall X,X' \in \cX \right \}.
}
Our main results provide general upper and lower bounds for the minimax risks \ef{minimax-fixed} and \ef{minimax-random} for these classes. However, before presenting these, it is worth noting the following, which is a trivial consequence of our lower bounds.
\prop{[Algebraic decay is impossible]
\label{prop:no-alg-conv}
Suppose that Assumption \ref{ass:mu} holds, as well as either Assumption \ref{ass:hilbert-noise} or Assumption \ref{ass:white-noise}. Let $1 \leq p < \infty$ and, in the latter case, $t > 0$. Then, for any $q > 0$ 
, we have
\bes{
\limsup_{m \rightarrow \infty}  \cM_{m}(\cF ; L^p_{\mu}) \cdot m^q  = + \infty,
}
where $\cF = \cF_{B,L}$ if Assumption holds \ref{ass:hilbert-noise} 
and $\cF = \cF^{t}_{B,L}$ otherwise.
}
This result demonstrates a \textit{curse of sample complexity}: The minimax error is necessarily subalgebraically decaying in $m$. Notice that this holds for \textit{any} measure $\mu$ satisfying Assumption \ref{ass:mu}---in particular, the curse of sample complexity arises regardless of the decay rate of the eigenvalues. 

With this result in mind, we next focus on establishing concrete upper and lower bounds for the minimax risk. As we shall see, these bounds depend critically on the decay of the eigenvalues $\{\lambda_i\}_{i\geq1}$. To obtain concrete guarantees we specialize to specific decay profiles of these eigenvalues (e.g., algebraic, exponential and double-exponential). 
A key result of this paper is a tight characterization for exponentially-decaying eigenvalues, which we now present.

\thm{[Tight characterization for exponentially-decaying $\lambda_i$; Hilbert-valued Gaussian noise]
\label{thm:main-exp-decay}
Suppose that Assumption \ref{ass:mu} and Assumption \ref{ass:hilbert-noise} hold. Let $1 \leq p < \infty$ and suppose that $\lambda_i = \exp(-\alphanew i^{\betanew})$ for some $\alphanew > 0$ and $\betanew \geq 1$. Then
\bes{
\cL_m(\cF_{B,L} ;  L^p_{\mu})  \asymp_{\alphanew,\betanew,\iota,p} \left(\log(m / \sigma^2) \right)^{\frac{\betanew}{\betanew+1}},\quad m \rightarrow \infty,
}
where $\cL_m(\cF_{B,L} ; L^p_{\mu}) = - \log(\cM_m(\cF_{B,L} ;  L^p_{\mu}))$ is as in \ef{neg-log-minimax}, with $\cM_m(\cF_{B,L} ;  L^p_{\mu})$ given by either  \ef{minimax-fixed} or \ef{minimax-random}.
}

\thm{[Tight characterization for exponentially-decaying $\lambda_i$; white noise]
\label{thm:main-exp-decay-white}
Suppose that Assumption \ref{ass:mu} and Assumption \ref{ass:white-noise} hold. Let $t > 0$, $1 \leq p < \infty$ and suppose that $w_i = i$ and $\lambda_i = \exp(-\alphanew i^{\betanew})$ for some $\alphanew > 0$ and $\betanew \geq 1$. Then
\bes{
\cL_m(\cF^t_{B,L} ; , L^p_{\mu})  \asymp_{\alphanew,\betanew,\iota,t,p} \left(\log(m / \sigma^2) \right)^{\frac{\betanew}{\betanew+1}},\quad m \rightarrow \infty.
}
where $\cL_m(\cF_{B,L} ; L^p_{\mu})$ is as in \ef{neg-log-minimax}, with $\cM_m(\cF_{B,L} ; , L^p_{\mu})$ given by either \ef{minimax-fixed} or \ef{minimax-random}.
}

These two results for the log-minimax risk imply that the corresponding minimax risk satisfies
\bes{
\cM_{m}(\cF_{B,L} ; L^p_{\mu}) \asymp C_1 \exp \left (- C_2 \left(\log(m / \sigma^2) \right)^{\frac{\betanew}{\betanew+1}} \right ),\quad m \rightarrow \infty,
}
where $C_1 = C_1(\alphanew,\betanew,\iota,\Upsilon,B,L,p)$ is a constant, potentially depending on all parameters, and $C_2 = C_2(\alphanew,\betanew,\iota,p)$. In particular, the minimax risk is \textit{subalgebraic}, but \textit{superlogarithmic}, regardless of the exponential decay rate $\betanew$ of the eigenvalues $\lambda_i$. Note that we do not specify the rate constant $C_2$ in this work. With additional work in tracking constants, our proofs would yield upper and lower bounds $C_2$. However, finding an exact expression for $C_2$ is an open problem.

\rem{
The constants in the lower and upper bounds implied by the $\asymp_{\alphanew,\betanew,\iota,p} $ in Theorem \ref{thm:main-exp-decay} symbol may only depend on only a subset of the parameters $\alphanew,\betanew,\iota,p$. For succinctness, we have not made this explicit. The specific parameters appearing in the lower and upper bounds can be readily inferred from the proofs of Theorem \ref{thm:main-exp-decay},\ref{thm:main-exp-decay-white}. Note that the same consideration applies to all subsequent results in this section.
}

These two results are consequences of general lower (Theorem \ref{t:main-lb}) and upper (Theorems \ref{t:main-ub-hilbert},\ref{t:main-ub-white}) bounds we establish later. These general bounds hold for arbitrary eigenvalues $\{\lambda_i\}_{i\geq1}$ and for both finite- and infinite-dimensional spaces $\cX$ and $\cY$. Like Theorems \ref{thm:main-exp-decay} and \ref{thm:main-exp-decay-white}, they also hold for arbitrary $\mu$ satisfying Assumption \ref{ass:mu}. Later, in Section \ref{sec:finite-dimensional-rates}, we show how classical minimax rates for bounded, Lipschitz functions on compact domains in $\bbR^d$ follow as simple consequences of these general upper and lower bounds.

\rem{[The case $0 < \betanew < 1$]
\label{rem:beta0-1}
Our tight characterization only holds for $\betanew \geq 1$. For $0 < \betanew < 1$, we have the non-tight bounds
\bes{
\left(\log(m / \sigma^2) \right)^{\frac{\betanew}{\betanew+1}} \lesssim_{\alphanew,\betanew,\iota,p} \cL_m(\cF_{B,L} ; L^p_{\mu}) \lesssim_{\alphanew,\betanew,\iota,p} \left(\log(m / \sigma^2) \right)^{\frac12}
}
in the case of Hilbert-valued noise. An analogous result holds for white noise. In particular, the upper bound for $\cL_m$ (i.e., lower bound for $\cM_m$) has exponent $1/2$ independent of $\betanew \in (0,1)$. We believe this is an artefact of the proof. See Remark \ref{rem:why-lower-bounds-nonsharp} for further discussion.
}

While our results tightly characterize the minimax risk in the case of exponentially decay, for slower-decaying eigenvalues we have non-matching upper and lower bounds. We now summarize these results in the case of algebraic decay. 

\thm{[Algebraically-decaying $\lambda_i$; Hilbert-valued Gaussian noise]
\label{thm:main-alg-decay}
Suppose that Assumption \ref{ass:mu} and Assumption \ref{ass:hilbert-noise} hold. Let $1 \leq p < \infty$ and suppose that $\lambda_i = i^{-\alphanew}$ for some $\alphanew > 1$. Then
\bes{
\frac{\alphanew-1}{2} \log \left ( \frac{\log(m/\sigma^2)}{\log(\log(m/\sigma^2))} \right ) + \cO_{\alphanew,\betanew,\iota,\Upsilon,B,L,p}(1) \leq \cL_m(\cF_{B,L} ; L^p_{\mu}) \lesssim_{\iota,p} \sqrt{\log(m/\sigma^2)},
}
as $m \rightarrow \infty$, where $\cL_m(\cF_{B,L} ; L^p_{\mu})$ is as in \ef{neg-log-minimax}, with $\cM_m(\cF_{B,L} ; , L^p_{\mu})$ given by either \ef{minimax-fixed} or \ef{minimax-random}.
}

\thm{[Algebraically-decaying $\lambda_i$; white noise]
\label{thm:main-alg-decay-white}
Suppose that Assumption \ref{ass:mu} and Assumption \ref{ass:white-noise} hold. Let $1 \leq p < \infty$ and suppose that $\lambda_i = i^{-\alphanew}$ for some $\alphanew > 1$. Then
\bes{
\frac{\alphanew-1}{2} \left ( \frac{\log(m/\sigma^2)}{\log(\log(m/\sigma^2))} \right ) + \cO_{\alphanew,\betanew,\iota,t,B,L,p}(1) \leq \cL_m(\cF^t_{B,L} ; L^p_{\mu}) \lesssim_{\iota,p} \sqrt{\log(m/\sigma^2)}.
}
where $\cL_m(\cF^t_{B,L} ; L^p_{\mu})$ is as in \ef{neg-log-minimax}, with $\cM_m(\cF^t_{B,L} ; , L^p_{\mu})$ given by either \ef{minimax-fixed} or \ef{minimax-random}.
}

These theorems imply the minimax risk satisfies, for either $\cF = \cF_{B,L}$ or $\cF = \cF^{t}_{B,L}$,
\bes{
C_1 \exp \left ( - C_2 \sqrt{\log(m/\sigma^2)}  \right ) \leq \cM_{m}(\cF ; L^p_{\mu}) \leq C_3 \left ( \frac{\log(m/\sigma^2)}{\log(\log(m/\sigma^2))} \right )^{-\frac{\alphanew-1}{2}}
}
for suitable constants $C_1,C_2,C_3 > 0$ depending on the parameters implied by the theorems. The lower bound is independent of algebraic decay rate $\alphanew$ of the eigenvalues---notably, it is the same rate as holds in the case of exponentially-decaying eigenvalues with $\betanew = 1$---and therefore we suspect it is nonsharp. In fact, we believe the minimax rate behaves like
\bes{
\left (\log(m/\sigma^2) \right )^{-\frac{\alphanew}{2}},\quad m \rightarrow \infty.
}
Namely, it is polylogarithmic in $m/\sigma^2$, albeit with a slightly better order than the corresponding bound given by Theorems \ref{thm:main-alg-decay},\ref{thm:main-alg-decay-white}. However, this remains an open problem.
The reason for the nontightness of the bounds in Theorems \ref{thm:main-alg-decay},\ref{thm:main-alg-decay-white} is the same as that discussed in Remark \ref{rem:beta0-1}. In particular, we refer to Remark \ref{rem:why-lower-bounds-nonsharp} for further discussion. 

To summarize, thus far in this section we have shown that minimax rates for Lipschitz operators are necessarily subalgebraic, and have giving a series of upper and lower bounds for different eigenvalues decay profiles. We remark in passing that, while impossible, rates that are close to algebraic can be obtained for sufficiently-fast decaying eigenvalues. In Section \ref{sec:double-exp-decay} we present a tight characterization for the log-minimax risk in the case of double-exponential decay, i.e., $\lambda_i = \exp(-\exp(\alphanew i))$. This corresponding minimax rate is algebraic for double exponentially-large regimes of $m/\sigma^2$.

The results in this section, which focus exclusively on Lipschitz operators, raise an important question: Does imposing higher regularity lead to faster, in particular algebraic, rates? In Section \ref{sec:double-exp-decay} we show that this is not the case. The same bounds hold for $C^{k,\alpha}$ (H\"older) operators for arbitrary $k \in \bbN_0$ and $0 < \alpha \leq 1$, with the only potential difference being in the constant. Thus imposing higher, but still finite regularity cannot overcome the curse of sample complexity.

\section{Lower bounds on the minimax risk} \label{sec:lower}

In this section, we establish the following lower bound on the minimax risk, from which the corresponding upper bounds in our main results will subsequently follow.

\thm{
\label{t:main-lb}
Suppose that Assumption \ref{ass:mu} holds and either Assumption \ref{ass:hilbert-noise} or \ref{ass:white-noise} hold. Let $\cF = \cF_{B,L}$ if Assumption \ref{ass:hilbert-noise} holds or $\cF = \cF^t_{B,L}$ for some $t > 0$ if Assumption \ref{ass:white-noise} holds. Then
\be{
\label{main-fixed-lb}
\cM_{m} (\cF ; L^p_{\mu}) \gtrsim \frac{L}{\sqrt{\sum^{d}_{j=1} 1/\lambda_j}}  \left ( \frac{\iota}{p+1} \right )^{d/p} \min \left \{ \frac{a}{8} , \left ( \frac{L^2 m}{c a^d \sigma^2 \sum^{d}_{j=1} 1/\lambda_j } \right )^{-\frac{1}{2+d}}  \right \} 
}
for all $d$ satisfying 
\be{
\label{d-cond-lb}
\frac{a L}{8 \sqrt{\sum^{d}_{i=1} 1/\lambda_i}} \leq B,
}
where $\cM_{m} (\cF ; L^p_{\mu})$ is given by either \ef{minimax-fixed} or \ef{minimax-random}. Here $c = c_{\Upsilon} > 0$ depends on $\Upsilon$ in the case of Assumption \ref{ass:hilbert-noise}, where $\Upsilon$ is the covariance operator of the noise, whereas $c > 0$ is universal in the case of Assumption \ref{ass:white-noise}.
}

\rem{
\label{rem:ben-lazy}
It is possible to remove the condition \ef{d-cond-lb}, at the expense of a more complicated lower bound \ef{main-fixed-lb}. However, this condition is not limiting in the settings considered in this paper. In particular, if $\sum^{\infty}_{i=1} \lambda_i = 1$ (i.e., the covariance operator has unit trace), then $\lambda_d \leq 1/d$ since the $\lambda_i$ are nonincreasing, and therefore this condition holds whenever $d \geq L / (8 B)$. In this work, when deriving lower bounds, we will always have $d \rightarrow \infty$ as $m \rightarrow \infty$, meaning that \ef{d-cond-lb} holds for all sufficiently large $m$.
}

\subsection{Key ingredients for the proof Theorem \ref{t:main-lb}}

Our lower bound follows the standard reduction of minimax estimation to a multi-hypothesis testing problem (see, e.g., Chapter~2 in~\cite{tsybakov2008introduction}). Fix a semidistance $d(\cdot,\cdot)$ on the model class $\mathcal{F}$. In this paper, we take
\be{
d(F,G) := \|F-G\|_{L^p_\mu(X;Y)}, \quad 1 \leq p < \infty.
}
For a finite subset of hypotheses $\{F_0,\dots,F_M\}\subset\mathcal{F}$, let $P_j$ denote the law of the observations under $F_j$. In the fixed-design setting, $P_j$ is the law of $(Y_1,\dots,Y_m)$. In the random-design setting, $P_j$ is the joint law of $\{(X_i,Y_i)\}_{i=1}^m$. If the hypotheses are well-separated, i.e.,
\be{
d(F_j,F_k)\ge 2s^*, \quad j\neq k,
}
and simultaneously the induced distributions are close in the Kullback--Leibler (KL) divergence, i.e.,
\be{
\frac{1}{M}\sum_{j=1}^M D(P_j\|P_0)\leq \alphanew\log M, \quad 0<\alphanew<1/8,
}
then no estimator can reliably distinguish the hypotheses. In which case, Fano's inequality provides a uniform lower bound of order $s^*$ for the minimax risk.

It remains to construct the well-separated set $\{F_0,\dots,F_M\} \subset \mathcal{F}$. We do this by combining (i) a geometric construction of many localized ``bump'' functions in the first $d$ eigencoordinates of $\mu$ with (ii) the Varshamov--Gilbert bound to extract a subset of $\{0,1\}^n$ with large Hamming separation, which translates into $L^p_\mu$-separation of the corresponding operators. We now recall the two standard tools (Fano's inequality and the Varshamov--Gilbert bound) that formalize this recipe.

\begin{theorem}[Fano's inequality]
\label{t:fano}
    Let $M \geq 2$ and let $\{F_0, \ldots, F_M\} \subset \cF$ be such that
    \begin{enumerate}
        \item $d(F_j, F_k) \geq 2 s^*$, where $d$ is a semidistance; and
        \item $\frac{1}{M} \sum_{j=1}^M D(P_j \| P_0) \leq \alphanew \log M$, with $0 < \kappa < 1/8$.
    \end{enumerate}
    Then
\be{
\label{Fano-prob}
        \inf_{\widehat{F}} \sup_{F \in \cF} \bbP \left [d(\widehat{F}, F) \geq s^* \right ] \geq \frac{\sqrt{M}}{1 + \sqrt{M}} \left ( 1 - 2\kappa - 2 \sqrt{\frac{\kappa}{\log M}} \right ) > 0
}
and, in particular,
\be{
\label{Fano-exp}
        \inf_{\widehat{F}} \sup_{F \in \cF} \bbE \left [ d(\widehat{F},F) \right ] \geq \frac{\sqrt{M}}{1 + \sqrt{M}} \left ( 1 - 2\kappa - 2 \sqrt{\frac{\kappa}{\log M}} \right ) s^*.
}
\end{theorem}

Note that \ef{Fano-exp} follows directly from \ef{Fano-prob} and Markov's inequality.
Next, let $\theta,\theta' \in \{0,1\}^n$. We define their \textit{Hamming distance} as
\bes{
H(\theta, \theta') = \sum_{j=1}^n \abs{\theta_j - \theta_j'} = \sum_{j=1}^n \abs{\theta_j - \theta_j'}^2.
}

\begin{lemma}[Varshamov--Gilbert bound]
\label{l:V-G}
    Let $n \geq 8$. Then there exists a subset $(\theta_0, \ldots, \theta_M)$ of $\{0,1\}^n$ such that 
\begin{enumerate}[label=(\roman*)]
        \item  $\theta_0 = (0, \ldots, 0)$;
        \item $M \geq 2^{n/8}$; and
        \item $H(\theta_j, \theta_k) \geq n/8$ for all $0 \leq j < k \leq M$.
    \end{enumerate}
\end{lemma}

\subsection{Proof of Theorem \ref{t:main-lb}}

To establish this result, we first require the following lemma.

\lem{
[Construction of a well-separated set]
\label{l:well-separated}
Let $0 < h / a \leq 1/8$ and $1 \leq p \leq \infty$, where $a$ is as in Assumption \ref{ass:mu}. Then there exists a collection of $L$-Lipschitz functionals $F_0,\ldots,F_M : \cX \rightarrow \bbR$ such that
\bes{
|F_j(X) | \leq \frac{L h}{\sqrt{\sum^{d}_{i=1} 1/\lambda_i}},\quad \forall X \in \cX,\ 0 \leq j \leq M,
}
and
\be{
\label{Fjk-lower}
\nm{F_j - F_{k}}_{L^p_{\mu}(\cX ; \bbR)} \gtrsim \left ( \frac{L}{\sqrt{\sum^{d}_{j=1} 1/\lambda_j}} \right ) \left ( \frac{\iota}{p+1} \right )^{d/p} h,\quad 0 \leq j < k \leq M,
}
where $M$ satisfies $\log(M) \geq (a/h)^{d} \log(2) / 8$. In the case $p = \infty$, the exponent $d/p$ in \ef{Fjk-lower} is interpreted as $0$.
}

\prf{
We divide the proof into several steps.

\pbk
\textit{Step 1: Bump function construction.} The first step is to construct a collection of bump functions. Let $d \in \bbN$ and define the operator $F : \cX \rightarrow \cY$ by
\bes{
F(X) =\prod^{d}_{i=1} \max \left \{ 1 - | \ip{X}{\phi_i}_{\cX} | / \sqrt{\lambda_i} , 0 \right \},
}
where we recall that the definitions of $\lambda_i$ and $\phi_i$ from Assumption \ref{ass:mu}.
This operator is Lipschitz with constant at most $\sqrt{\sum^{d}_{i=1} 1/\lambda_i}$.
Now $n \in \bbN$ and $C_1,\ldots,C_n \in \cX$ to be chosen later. Then, for $\theta \in \{0,1\}^n$, define the bump functions
\bes{
F_{\theta}(\cdot) = \frac{L h}{\sqrt{\sum^{d}_{i=1} 1/\lambda_i}} \sum^{n}_{i=1} \theta_i F \left ( \frac{\cdot - C_i }{h} \right ).
}
Observe that
\bes{
\mathrm{supp}\left ( F \left ( \frac{\cdot - C_i }{h} \right ) \right ) = \cE_h(C_i) : = \left \{ X \in \cX : | \ip{X - C_i}{\phi_j}_{\cX} |^2 \leq \lambda_j h^2,\ j = 1,\ldots,d \right \}.
}
We now make the following assumption: $h$ and $C_1,\ldots,C_n$ are chosen so that the sets $\cE_h(C_i)$ are disjoint and contained in the set
\bes{
\cE = \left \{ X \in \cX :  | \ip{X}{\phi_j}_{\cX} |^2 \leq a^2 \lambda_j,\ j = 1,\ldots,d \right \},
}
where $a>0$ is as in Assumption \ref{ass:mu}.
In that case, we immediately deduce that $F_{\theta}$ is $L$-Lipschitz. Moreover, since $|F(X) | \leq 1$,  
$\forall X \in \cX$, we have
\bes{
|F_{\theta}(X) | \leq \frac{L h}{\sqrt{\sum^{d}_{i=1} 1/\lambda_i}},\quad \forall X \in \cX,
}
as required.

\pbk
\textit{Step 2: Estimating the $L^2_{\mu}$-distance of two bump functions.} Let $\theta,\theta' \in \{0,1\}^n$. Then, by disjointness of the supports, we have, for any $1 \leq p < \infty$,
\be{
\label{L2-theta-split}
\begin{split}
\nm{F_{\theta} - F_{\theta'}}^p_{L^p_{\mu}(\cX ; \bbR)} &= \left ( \frac{L h}{\sqrt{\sum^{d}_{j=1} 1/\lambda_j}} \right )^p \sum^{n}_{i=1} |\theta_i - \theta'_i |^p \int_{\cX} \left | F \left ( \frac{X- C_i }{h} \right ) \right |^p \D \mu(X) 
\\
& = :  \left ( \frac{L h}{\sqrt{\sum^{d}_{j=1} 1/\lambda_j}} \right )^p \sum^{n}_{i=1} |\theta_i - \theta'_i| I_{i,p}. 
\end{split}
}
Here, in the second step, we also used the fact that $\theta , \theta' \in \{0,1\}^n$, which implies that $|\theta_i - \theta'_i|^p = | \theta_i - \theta'_i|$.
Recall that $\{ \phi_i \}_{i \in \bbN}$ is an orthonormal basis of $\cX$. Notice that the function $F(X) = f \circ T(X)$, where
\bes{
T : \cX \rightarrow \bbR^d,\ X \mapsto (\ip{X}{\phi_i}_{\cX} / \sqrt{\lambda_i} )^{d}_{i=1}
}
and
\bes{
f : \bbR^d \rightarrow \bbR,\ x \mapsto \prod^{d}_{i=1} \max \{ 1 - |x_i| , 0 \}.
}
Let $c_i = T(C_i)$. Then
\bes{
I_{i,p} = \int_{\bbR^d} \left | f \left ( \frac{x - c_i}{h} \right ) \right |^p \D \mu \sharp T^{-1}(x).
}
Assumption \ref{ass:mu} gives that $\mu \sharp T^{-1} = \bigotimes^{d}_{i=1} \xi_i$.
Therefore
\bes{
I_{i,p} = \int_{\bbR^d} \left | f \left ( \frac{x - c_i}{h} \right ) \right |^p \prod^{d}_{i=1} \nu_i(x) \D x.
}
where $\nu_i$ is the density of the random variable $\xi_i$.
Notice that 
\bes{
\mathrm{supp} \left ( f \left ( \frac{x-c_i}{h} \right ) \right ) = D_h(c_i) = \left \{ x \in \bbR^d : (x_j - (c_i)_j)^2 \leq h^2,\ j = 1,\ldots,d \right \}
}
and also, by assumption, that
\bes{
D_h(c_i) \subseteq D : = \left \{ x \in \bbR^d : |x_j| \leq a,\ j = 1,\ldots,d \right \}.
}
We now use this and Assumption \ref{ass:mu} to write
\bes{
I_{i,p} \geq b^d \int_{D} \left | f \left ( \frac{x - c_i}{h} \right ) \right |^p  \D x = (bh)^d \left ( \int^{1}_{-1} (1-|x|)^p \D x \right )^d = h^d (2b/(p+1))^d.
}
Substituting this into \ef{L2-theta-split}, we deduce that
\bes{
\nm{F_{\theta} - F_{\theta'}}^p_{L^p_{\mu}(\cX ; \bbR)} \geq  \left ( \frac{L }{\sqrt{\sum^{d}_{j=1} 1/\lambda_j}} \right )^p  \left(\frac{2b}{p+1} \right )^d h^{d+p} H(\theta,\theta'). 
}
where $H$ is the Hamming distance, i.e.,
\bes{
\nm{F_{\theta} - F_{\theta'}}_{L^p_{\mu}(\cX ; \bbR)} \geq \left ( \frac{L }{\sqrt{\sum^{d}_{j=1} 1/\lambda_j}} \right ) \left(\frac{2b}{p+1} \right )^{d/p}
h^{d/p+1} (H(\theta,\theta'))^{1/p}. 
}
This holds for $1 \leq p < \infty$. For $p = \infty$, we have
\bes{
\nm{F_{\theta} - F_{\theta'}}_{L^{\infty}_{\mu}(\cX ; \bbR)} = \left ( \frac{L}{\sqrt{\sum^{d}_{i=1} 1/\lambda_i} } \right ) h\max_{i=1,\ldots,d} | \theta_i - \theta'_i | =  \left ( \frac{L}{\sqrt{\sum^{d}_{i=1} 1/\lambda_i} } \right ) h ,\quad \forall \theta' \neq \theta. 
}

\pbk
\textit{Step 3: Extracting a subset.} 
Now suppose that $n \geq 8$. Then Lemma \ref{l:V-G} implies there exists a collection $\theta_0,\ldots,\theta_M \in \{0,1\}^n$ such that $M \geq 2^{n/8}$ and $H(\theta_j,\theta_k) \geq n/8$ for $j \neq k$. Let $F_j = F_{\theta_j}$, $j = 0,\ldots,M$. Then we deduce that
\be{
\label{Fjk-dist-n-h}
\nm{F_j - F_{k}}_{L^p_{\mu}(\cX ; \bbR)} \gtrsim \left ( \frac{L}{\sqrt{\sum^{d}_{j=1} 1/\lambda_j}} \right ) \left ( \frac{2 b}{p+1} \right )^{d/p} h^{d/p+1} n^{1/p}.
} 
This also holds when $p = \infty$, interpreting $1/p$ as $0$.

\pbk
\textit{Step 4: Choosing $n$ and $C_1,\ldots,C_n$.} We require that $n \geq 8$ and the sets $\cE_h(C_i)$ are disjoint and contained in $\cE$. The latter assumption is equivalent to assuming that the sets $D_h(c_i)$ are disjoint and contained in $D$. This is equivalent to the question of packing the $\infty$-norm ball $D = \{ x \in \bbR^d : \nm{x}_{\infty} \leq a\}$ in $\bbR^d$ with $\infty$-norm balls of radius $h$. Hence, we can choose $n$ as the packing number of $D$ with radius $h$. This packing number is bounded below by $(a/h)^d$. Since $h \leq a/8$ by assumption, we have $n \geq (a/h)^d \geq a/h \geq 8$, as required. Substituting this into \ef{Fjk-dist-n-h} and recalling that $M \geq 2^{n/8}$ now yields the result.
}

\prf{[Proof of Theorem \ref{t:main-lb}] 
Let $0 < h/a \leq 1/8$. We combine the previous lemma with Fano's inequality (Theorem \ref{t:fano}). Consider the functionals $F_0,\ldots,F_M$ from the previous lemma. We use these to construct operators by multiplying by a suitable element $Y \in \cY$: Namely, $G_j(X) = Y F_j(X)$ for $X \in \cX$ and $j = 0,\ldots,M$. The choice of $Y$ is dictated by the noise model as follows:
\begin{enumerate}[label=(\roman*)] 
\item For Hilbert-valued Gaussian noise (Assumption \ref{ass:hilbert-noise}) let $Y = Y_1$, where $Y_1$ is the eigenvector of the covariance operator $\Upsilon$ corresponding to its largest eigenvalue $\upsilon_1$. 
\item For white noise (Assumption \ref{ass:white-noise}), we set $Y = \psi_1$, where $\psi_1$ is the first element of the orthonormal basis that defines the scale space $\cY^t$ in \ef{Yt-def}.
\end{enumerate}
Observe that $\nm{Y}_{\cY} = 1$ in both cases. Recall that $h/a \leq 1/8$, by assumption. Therefore, the $F_j$ satisfy
\bes{
|F_j(X)| \leq \frac{L h}{\sqrt{\sum^{d}_{j=1} 1/\lambda_j}} \leq \frac{a L}{8 \sqrt{\sum^{d}_{j=1} 1/\lambda_j}} \leq B,
}
where the final inequality is due to \ef{d-cond-lb}.
Hence, in (i) we have $G_j \in \cF_{B,L}$, $\forall j$, and in (ii) we have $G_j \in \cF^{t}_{B,L}$, $\forall j$. Let $d(\cdot,\cdot)$ denote the $L^p_{\mu}$-distance, set $\kappa = 1/16$ (this choice is arbitrary) and set
\bes{
s^* =  \frac{c_0}{2}\frac{L}{\sqrt{\sum^{d}_{j=1} 1/\lambda_j}} \left ( \frac{\iota}{p+1} \right )^{d/p} h,
}
where $c_0>0$ is the universal constant in \ef{Fjk-lower}. Notice that $\log(M) \geq (a/h)^d \log(2) / 8 \geq \log(2)$, since $0 < h/a \leq 1/8$.

Consider the fixed designs case. In case (i), we have 
\bes{
D(P_j \| P_0 ) = \frac{1}{2\sigma^2} \sum_{i=1}^m \nm{\Upsilon^{-1/2} G_j(X_i)}^2_{\cY} \leq \frac{m}{2\sigma^2 \upsilon_1} \max_{X \in \cX} | F_j(X)|^2 \leq c \frac{L^2 h^2 m}{\sigma^2 \sum^{d}_{j=1} 1/\lambda_j },
}
where $c = c_{\Upsilon} > 0$ depends on $\Upsilon$ (in fact, its maximum eigenvalue) only. In case (ii), we have
\bes{
D(P_j \| P_0 ) = \frac{1}{2\sigma^2} \sum_{i=1}^m \nm{G_j(X_i)}^2_{\cY} \leq \frac{m}{2 \sigma^2}\max_{X \in \cX} | F_j(X)|^2 \leq c \frac{L^2 h^2 m}{\sigma^2 \sum^{d}_{j=1} 1/\lambda_j },
}
where $c = 1/2 > 0$ is a universal constant. The arguments for random designs are similar. Indeed, we have
\bes{
D(P_j \| P_0) \leq c \frac{m}{\sigma^2} \nm{F_j}^2_{L^2_{\nu}(\cX ; \bbR)} \leq c \frac{m}{\sigma^2} \max_{X \in \cX} | F_j(X)|^2 \leq c \frac{L^2 h^2 m}{\sigma^2 \sum^{d}_{j=1} 1/\lambda_j },
}
for the same constant $c >0$ as before. 
Therefore, in all settings, we deduce that
\bes{
\frac1M \sum^{M}_{j=1} D(P_j \| P_0 ) \leq c \frac{L^2 h^2 m}{\sigma^2  \sum^{d}_{j=1} 1/\lambda_j }.
}
Applying Theorem \ref{t:fano} and using the fact that $\log(M) \geq (a/h)^d \log(2)/8$, we deduce that
\bes{
\cM_m(\cF;L^p_{\mu}) \gtrsim \frac{L}{\sqrt{\sum^{d}_{j=1} 1/\lambda_j}} \left ( \frac{\iota}{p+1} \right )^{d/p} h ,
}
provided
\bes{
\frac{L^2 h^2 m}{\sigma^2  \sum^{d}_{j=1} 1/\lambda_j} \leq c (a/h)^{d},
}
for a potentially different constant $c$, where $c = c_{\Upsilon} > 0$ in case (i) and $c >0$ is universal in case (ii).
Rearranging, we now set
\bes{
h = \min \left \{ \frac{a}{8} , \left ( \frac{L^2 m}{ c a^d \sigma^2 \sum^{d}_{j=1} 1/\lambda_j } \right )^{-\frac{1}{2+d}} \right \}.
}
Substituting this into the previous expression now gives the result.
}

\section{Upper bounds on the minimax risk} \label{sec:upper}

We now establish upper bounds on the minimax risks. For these, we split into two cases corresponding to the different noise distributions (Assumption \ref{ass:hilbert-noise} or \ref{ass:white-noise}).

\thm{
\label{t:main-ub-hilbert}
Suppose that Assumptions \ref{ass:mu} and \ref{ass:hilbert-noise} hold. Let $m \in \bbN$. Then
there is a fixed design $\cD = \{X_1,\ldots,X_m\}$ (or, in the random design case, a random design $\nu$) and an estimator satisfying
\bes{
\bbE\nm{F - \widehat{F}}_{L^p_{\mu}(\cX ; \cY)} \lesssim_p  \left ( \frac{m}{\sigma^2 \sqrt{\lambda_1\cdots\lambda_d} } \right )^{-\frac{2}{(p+2)d+4} } (B^{2p} L^4 d^4 )^{\frac{d}{2((p+2)d+4)}} +L \sqrt{\sum_{j>d} \lambda_j}
}
for all $F \in \cF_{B,L}$ and all $d \in \bbN$ for which
\bes{
 \frac{m}{\sigma^2}   \geq \frac{(\lambda_1\cdots\lambda_d)^{\frac{1}{2}}}{\left ( B^{2p} L^4 d^4 \right )^{\frac{1}{p+2}} (\lambda_d)^{\frac{(p+2) d+4}{2(p+2)}} }. 
}
}

The white noise case is similar, except for an additional term in the error bound depending on the smoothness space \ef{Yt-def}.

\thm{
\label{t:main-ub-white}
Suppose that Assumptions \ref{ass:mu} and \ref{ass:white-noise} hold. Then for each $m \in \bbN$ there is a fixed design $\cD = \{X_1,\ldots,X_m\}$  (or, in the random design case, a random design $\nu$) such that, for each $r \in \bbN$, there there is an estimator satisfying
\bes{
\bbE \nm{F - \widehat{F}}_{L^p_{\mu}(\cX ; \cY)} \lesssim_p B w^{-t}_r + \left ( \frac{m}{r \sigma^2 \sqrt{\lambda_1\cdots\lambda_d} } \right )^{-\frac{2}{(p+2)d+4} } (B^{2p} L^4 d^4 )^{\frac{d}{2((p+2)d+4)}} +L \sqrt{\sum_{j>d} \lambda_j}
}
for all $F \in \cF_{B,L}$ and all $d \in \bbN$ for which
\bes{
 \frac{m}{\sigma^2}   \geq \frac{r (\lambda_1\cdots\lambda_d)^{\frac{1}{2}}}{\left ( B^{2p} L^4 d^4 \right )^{\frac{1}{p+2}} (\lambda_d)^{\frac{(p+2) d+4}{2(p+2)}} }. 
}
}

\subsection{Construction of the estimator}

A standard approach for establishing upper bounds for Lipschitz functions in finite dimensions and for compactly-supported measures $\mu$ involves constructing a \textit{histogram} estimator. We follow a similar approach, extending this estimator to infinite-dimensional Lipschitz operators and measures with potentially noncompact supports. 

Let $d \in \bbN$, $n_1,\ldots,n_d \in \bbN$ and $R > 0$ be parameters. 
Let $\lambda_i$ and $\phi_i$, $i \in \bbN$, be as in Assumption \ref{ass:mu} and define the set
\bes{
D = \left \{ X \in \cX : | \ip{X}{\phi_i}_{\cX} | \leq \sqrt{R \lambda_i},\ \forall i =1,\ldots,d \right \} \subset \cX.
}
Now subdivide each interval $[-\sqrt{R \lambda_i} , \sqrt{R \lambda_i } ]$ into $n_i$ equally-spaced intervals $I_{i,j}$, $j = 1,\ldots,n_i$, and consider the cells
\bes{
C = C_{j_1,\ldots,j_d} = \left \{ X \in \cX : \ip{X}{\phi_i}_{\cX} \in I_{i,j_i},\ \forall i = 1,\ldots,d \right \},
}
Let $\cC$ denote the collection of cells. Note that
\bes{
|\cC| = n_1 \cdots n_d = : n.
}
Given $F : \cX \rightarrow \cY$, in the case of Hilbert-valued Gaussian noise, we define the estimator
\be{
\label{Fhat-hilbert-noise}
\widehat{F} = \sum_{C \in \cC} Y_{C} \bbI_{C},
}
where $\bbI_C$ denotes the indicator function of the set $C$ and
\bes{
Y_{C} = \frac{1}{N_C} \sum_{i : X_i \in C} Y_i\quad \text{and}\quad N_C = | \{ i : X_i \in C \} |.
}
For completeness, if $N_C = 0$ we define $Y_C = 0$.
In the case of white noise, we modify the estimator as follows. Let $r \in \bbN$ and define the operator
\bes{
\cS_r : \cY \rightarrow \cY,\ Y \mapsto \sum^{r}_{i=1} \ip{Y}{\psi_i}_{\cY} \psi_i,
}
where $\{ \psi_i \}^{\infty}_{i=1}$
is the orthonormal basis of $\cY$ that defines the spaces $\cY^t$ in \ef{Yt-def}.
Then we set
\be{
\label{Fhat-white-noise}
\widehat{F} = \sum_{C \in \cC} Y_{C,r} \bbI_{C},\quad \text{where }Y_{C,r} = \cS_r(Y_C).
}

\subsection{Error bounds for the estimator}

\lem{
\label{lem:fixed-design-histogram-bd}
Consider the fixed design setting, where the $X_1,\ldots,X_m$ are such that $N_C = m/n$ (assumed to be an integer) and $\ip{X_i}{\phi_j}_{\cX} =0$, $\forall j > d$. In the Hilbert-valued Gaussian noise setting, the estimator \ef{Fhat-hilbert-noise} satisfies, for any $F \in \cF_{B,L}$,
\bes{
\bbE\nm{F - \widehat{F}}_{L^p_{\mu}(\cX ; \cY)} \lesssim_p  \sqrt{\frac{n}{m}} \sigma +  L \sqrt{R} \sqrt{\sum^{d}_{j=1} \frac{\lambda_j}{n^2_j}} + L \sqrt{\sum_{j > d} \lambda_j} + B \mu(D^c)^{\frac1p}.
}
In the white noise setting, the estimator \ef{Fhat-white-noise} satisfies, for any $t \geq 0$ and $F \in \cF^{t}_{B,L}$,
\bes{
\bbE \nm{F - \widehat{F}}_{L^p_{\mu}(\cX ; \cY)} \lesssim_p B w^{-t}_{r} + \sqrt{\frac{r n}{m}} \sigma + L \sqrt{R} \sqrt{\sum^{d}_{j=1} \frac{\lambda_j}{n^2_j}} + L \sqrt{\sum_{j > d} \lambda_j} + B \mu(D^c)^{\frac1p}.
}
}
\prf{
Since $\cS_r \rightarrow \cI$ converges strongly to the identity operator on $\cY$, we may consider the estimator \ef{Fhat-hilbert-noise} as a special case of \ef{Fhat-white-noise} with $r = \infty$ and $\{ \psi_i \}_{i}$ being the orthonormal basis of eigenvectors of the covariance of Hilbert-valued Gaussian
noise distribution. Assuming $r \in \bbN \cup \{ \infty \}$, we first write
\bes{
\widetilde{F} = \sum_{C \in \mathcal{C}} \widetilde{Y}_{C,r} \bbI_C,\quad \widetilde{Y}_{C,r} = \frac{1}{N_C} \sum_{i : X_i \in C} \cS_r(F(X_i)),
}
and then write
\be{
\label{ub-err-split}
\nm{F - \widehat{F}}_{L^p_{\mu}(\cX ; \cY)} \leq \nm{F - \cS_r \circ F}_{L^p_{\mu}(\cX ; \cY)} +  \nm{\cS_r \circ F - \widetilde{F}}_{L^p_{\mu}(\cX ; \cY)} + \nm{\widetilde{F} - \widehat{F}}_{L^p_{\mu}(\cX ; \cY)}.
}
Consider the first term. We have
\eas{
\nm{F(X) - \cS_r \circ F(X)}^2_{\cY} &= \sum_{i > r} | \ip{F(X)}{\psi_i}_{\cY} |^2
\\
& \leq w^{-2t}_r \sum_{i > r} w^{2t}_r | \ip{F(X)}{\psi_i}_{\cY} |^2  
\\
& \leq w^{-2t}_r \nm{F(X)}^2_{\cY^t} 
\\
& \leq B^2 w^{-2t}_r .
}
Therefore
\be{
\label{ub-err-term1}
\nm{F - \cS_r \circ F}_{L^p_{\mu}(\cX ; \cY)} \leq B w^{-t}_{r} .
}
Consider the third term of \ef{ub-err-split}. We have 
\eas{
\nm{\widetilde{F} - \widehat{F}}^p_{L^p_{\mu}(\cX ; \cY)} & = \sum_{C \in \cC} \mu(C) \nm{Y_{C,r} - \widetilde{Y}_{C,r} }^p_{\cY}  =  \sum_{C \in \cC} \mu(C) \frac{\sigma^p}{N^p_C} \nms{ \sum_{i : X_i \in C} \cS_r(E_i)}^p_{\cY} .
}
Consider the Hilbert-valued noise case. Let $\{ \psi_j \}_{j \in \bbN} \subset \cY$ be the orthonormal basis of eigenvectors of the covariance operator $\Upsilon$ of the noise and $\upsilon_1 \geq \upsilon_2 \geq \cdots \geq 0$ be the eigenvalues. Let
\bes{
I_C = \bbE \nms{ \sum_{i : X_i \in C} E_i }^p_{\cY} = \bbE \left ( \sum^{\infty}_{j=1} \left | \sum_{i : X_i \in C} \ip{E_i}{\psi_j}_{\cY} \right |^2 \right )^{\frac{p}{2}} .
}
Now observe that $\ip{E_i}{\psi_j} \sim \cN(0,\upsilon_j)$ 
and therefore we can write
\bes{
\sum_{i : X_i \in C}\ip{E_i}{\psi_j}_{\cY}  = \sigma \sqrt{\upsilon_j} Z_j,\quad \text{where } Z_j = \sum_{i : X_i \in C} n_{ij}.
}
and $n_{ij} \sim_{\mathrm{i.i.d.}} \cN(0,1)$. We now apply H\"older's inequality to get
\eas{
I_C &= \bbE \left [ \sum^{\infty}_{j=1} \upsilon_j |Z_j|^2 \right ]^{\frac{p}{2}} 
\\
& \leq  \bbE \left [ \left ( \sum^{\infty}_{j=1} \upsilon^{(1-1/p)q}_j \right )^{\frac1q} \left ( \sum^{\infty}_{j=1} \upsilon_j |Z_j|^{2p} \right )^{\frac1p} \right ]^{\frac{p}{2}} 
\\
&=  \mathrm{Tr}(\Upsilon)^{\frac{p-1}{2}} \bbE \left ( \sum^{\infty}_{j=1} \upsilon_j |Z_j|^{2p} \right )^{\frac12}
\\
& \leq  \mathrm{Tr}(\Upsilon)^{\frac{p-1}{2}} \left ( \sum^{\infty}_{j=1} \upsilon_j \bbE |Z_j|^{2p} \right )^{\frac12},
}
where $1 < q \leq \infty$ satisfies $1/p+1/q=1$. Now $Z_j$ is a sum of $N_C$ i.i.d.\ $\cN(0,1)$ random variables. Therefore, it is subgaussian with parameters $\betanew = 2$ and $\kappa = 1/(2 N_C)$. By a standard moment bound, we have
\bes{
[\bbE ( |Z_j |^{2p} ) ]^{\frac{1}{2p}} \lesssim_p \sqrt{N_c}.
}
We deduce that
\bes{
I_C \lesssim_p  \mathrm{Tr}(\Upsilon)^{\frac{p}{2}} N^{\frac{p}{2}}_C .
}
Therefore
\be{
\label{ub-err-term2-partial}
\bbE\nm{\widetilde{F} - \widehat{F}}_{L^p_{\mu}(\cX ; \cY)} \lesssim_p \left ( \sigma^{p} \mathrm{Tr}(\Upsilon)^{\frac{p}{2}} \sum_{C \in \cC} \frac{\mu(C)}{N^{\frac{p}{2}}_C} \right )^{\frac1p} \leq \sigma \sqrt{\frac{\mathrm{Tr}(\Upsilon) n}{m}},
}
where we recall that $N_C = m/n$. In the Hilbert-valued noise case, since $\mathrm{Tr}(\Upsilon) = 1$, this gives 
\be{
\label{ub-err-term2a}
\bbE \nm{\widetilde{F} - \widehat{F}}_{L^p_{\mu}(\cX ; \cY)} \lesssim_p  \sigma \sqrt{\frac{n}{m}}.
}
However, this analysis also applies to the white noise case, upon setting $\upsilon_1 = \cdots = \upsilon_r = 1$ and $\upsilon_i =0$ otherwise. We deduce that 
\be{
\label{ub-err-term2b}
\bbE \nm{\widetilde{F} - \widehat{F}}_{L^p_{\mu}(\cX ; \cY)} \lesssim_p  \sigma \sqrt{\frac{r n}{m}}
}
in the white noise case.

Now consider the second term in \ef{ub-err-split}. We first write
\bes{
\nm{\cS_r \circ F - \widetilde{F} }^p_{L^p_{\mu}(\cX ; \cY)} = \int_{D} \nm{\cS_r \circ F(X) - \widetilde{F}(X)}^p_{\cY} \D \mu(X) +  \int_{D^c} \nm{\cS_r \circ F(X)}^p_{\cY} \D \mu(X).
}
For the second term, since $\cS_r$ is a projection, we have
\bes{
\nm{\cS_r \circ F(X)}_{\cY} \leq \nm{F(X)}_{\cY} \leq \nm{F(X)}_{\cY^t} \leq B
}
and therefore 
\bes{
\int_{D^c} \nm{\cS_r \circ F(X)}^p_{\cY} \D \mu(X) \leq \mu(D^c) B^p .
}
Next, let $X \in C$ for some $C \in \cC$. Then, since $\cS_r$ is a projection and $F$ is $L$-Lipschitz, we have
\eas{
\nm{\cS_r \circ F(X) - \widetilde{F}(X)}_{\cY} \leq \frac{1}{N_C} \sum_{i : X_i \in C} \nm{\cS_r(F(X)) - \cS_r(F(X_i))}_{\cY} \leq L \max_{i : X_i \in C} \nm{X - X_i}_{\cX}.
}
Now, by construction, 
\bes{
\nm{X - X_i}^2_{\cX} = \sum^{d}_{j=1} | \ip{X - X_i}{\phi_j}_{\cX} |^2 + \sum_{j > d} | \ip{X}{\phi_j}_{\cX} |^2 \leq 4 R \sum^{d}_{j=1} \frac{\lambda_j}{n^2_j} + \sum_{j > d} | \ip{X}{\phi_j}_{\cX} |^2 .
}
Therefore
\eas{
\nm{\cS_r \circ F - \widetilde{F}}^p_{L^p_{\mu}(\cX ; \cY)}  \lesssim_p \left ( L^2 R \sum^{d}_{j=1} \frac{\lambda_j}{n^2_j} \right )^{p/2} + L^p  \int_{\cX} \left(\sum_{j > d}| \ip{X}{\phi_j}_{\cX} |^2 \right )^{p/2} \D \mu(X) + \mu(D^c) B^p .
}
Consider the second term, which we denote as $I$ for convenience. By assumption, $\ip{X}{\phi_j}_{\cX} = \sqrt{\lambda_j} \xi_j$. Hence we may write this as
\bes{
I = \bbE \left ( \sum_{j > d} \lambda_j | \xi_j |^2 \right )^{p/2},
}
where the expectation is with respect to the $\xi_j$. Suppose first that $p \geq 2$. Then, by Minkoswki's inequality and the moment assumption on the $\xi_j$, we have
\bes{
I^{\frac{2}{p}} \leq \sum_{j > d} \left [ \bbE (\lambda_j | \xi_j |^2 )^{\frac{p}{2}} \right ]^{\frac{2}{p}} = \sum_{j > d} \lambda_j \left [ \bbE | \xi_j |^{p} \right ]^{\frac{2}{p}} \lesssim_{\xi,p} \sum_{j > d} \lambda_j
}
and therefore
\bes{
I^{\frac1p} \lesssim_{\xi,p} \sqrt{\sum_{j > d} \lambda_j}.
}
Now suppose that $1 \leq p \leq 2$. Then
\bes{
I \leq \left [ \bbE \left ( \sum_{j > d} \lambda_j |\xi_j|^2 \right ) \right ]^{\frac{p}{2}} = \left ( \sum_{j > d} \lambda_j \bbE |\xi_j|^2 \right )^{\frac{p}{2}} \lesssim_{\xi,p} \left ( \sum_{j > d} \lambda_j \right )^{\frac{p}{2}},
}
and therefore, once more,
\bes{
I^{\frac1p} \lesssim_{\xi,p} \sqrt{\sum_{j > d} \lambda_j}.
}
Substituting this into the previous expression, we see that
\be{
\label{ub-err-term3}
\nm{\cS_r \circ F - \widetilde{F}}_{L^p_{\mu}(\cX ; \cY)}  \lesssim_{\xi,p} L \sqrt{R}\sqrt{\sum^{d}_{j=1} \frac{\lambda_j}{n^2_j}} + L \sqrt{\sum_{j > d} \lambda_j } + B (\mu(D^c))^{\frac1p}.
}
Substituting this, \ef{ub-err-term1} and either \ef{ub-err-term2a} or \ef{ub-err-term2b} into \ef{ub-err-split} now gives the result.
}

We now consider the random design setting. The following lemma shows the existence of a random design $\nu$ that achieves the same error bound.

\lem{
\label{lem:random-design-histogram-bd}
Consider the probability distribution $\nu$ on $\cX$ given by the law of the random variable
\bes{
X = \sum^{d}_{i=1} \sqrt{R \lambda_i} \zeta_i \phi_i,\quad \text{where }\zeta_i \sim_{\mathrm{i.i.d.}} \cU([-1,1]).
}
Let $X_1,\ldots,X_m \sim_{\mathrm{i.i.d.}} \nu$.
Then the same conclusions as Lemma \ref{lem:fixed-design-histogram-bd} hold, where the expectation is now also taken with respect to design points $X_i$ in addition to the noise.
}
\prf{
We modify the proof of the previous lemma. 
Consider the three terms on the right-hand side of \ef{ub-err-split}. It follows from the derivations leading to \ef{ub-err-term1} and \ef{ub-err-term3} that these bounds for first and second terms of \ef{ub-err-split} are unaffected choice of design points. Therefore, we only need to show that the bounds \ef{ub-err-term2a},\ef{ub-err-term2b} for the third term remain valid for this choice of random design. Following the same arguments, we see that
\bes{
\bbE \nm{\widetilde{F} - \widehat{F}}^p_{L^2_{\mu}(\cX ; \cY)} \lesssim \sigma^p \mathrm{Tr}(\Upsilon)^{\frac{p}{2}} \sum_{C \in \cC} \mu(C) \bbE \left ( \frac{1}{N^{\frac{p}{2}}_C} \right ),
}
where we observe that $N_C$ is now the random variable given by $N_C = | \{ i : X_i \in C \} |$ whenever at least one $X_i$ belongs to $C$ and $+ \infty$ otherwise. Let 
\bes{
p_C = \bbP_{X \sim \nu}(X \in C) = \frac{1}{n},
}
where the second equality follows from the definition of $\nu$. Then $N_C \sim \mathrm{Bin}(m,p_C)$.  A standard Chernoff bound gives that
\bes{
\bbP\left( N_C \leq \frac{m }{2 n} \right ) \leq \exp \left (-\frac{m }{8 n} \right ).
}
Therefore, we may write
\eas{
\bbE \left ( \frac{1}{N^{p/2}_C} \right ) &= \bbE \left ( \frac{1}{N^{p/2}_C} \Bigg | N_C \leq \frac{m}{2 n} \right ) \bbP \left ( N_C \leq \frac{m}{2n} \right ) + \bbE \left (  \frac{1}{N^{p/2}_C} \Bigg | N_C > \frac{m}{2n} \right )  \bbP \left ( N_C > \frac{m}{2n} \right )
\\
& \leq \exp \left ( - \frac{m}{8n} \right ) + \left ( \frac{m}{2n} \right )^{-\frac{p}{2}}
\\
& \lesssim_p \left (\frac{m}{n} \right )^{-\frac{p}{2}},
}
where in the last step we used the fact that the function $x^{p/2} \exp(-x)$ is bounded on $[0,\infty)$. This implies that
\bes{
\bbE \nm{\widetilde{F} - \widehat{F}}_{L^2_{\mu}(\cX ; \cY)} \lesssim_p \sigma \sqrt{\frac{\mathrm{Tr}(\Upsilon) n}{m}}.
}
In other words, \ef{ub-err-term2-partial} remains valid in this setting. We immediately deduce that \ef{ub-err-term2a},\ef{ub-err-term2b} also hold. This completes the proof.
}

\subsection{Proofs of Theorems \ref{t:main-ub-hilbert} and \ref{t:main-ub-white}}

\prf{[Proof of Theorem \ref{t:main-ub-hilbert}]
We use Lemma \ref{lem:fixed-design-histogram-bd} (for fixed designs) or Lemma \ref{lem:random-design-histogram-bd} (for random designs). First, we estimate $\mu(D^c)$. By Assumption \ref{ass:mu}(i)--(ii), the union bound and Markov's inequality, we have
\bes{
\mu(D^c) \leq \sum^{d}_{i=1} \bbP_{x \sim \xi_i}[|x|^2 \geq R \lambda_i ] \leq \sum^{d}_{i=1} \frac{\bbE_{x \sim \xi_i} |x|^2}{R \lambda_i} = \frac{d}{R}.
}
We now optimize with respect to $R$, by setting
\bes{
R = \left ( \frac{B d^{\frac1p}}{L \sqrt{\sum^{d}_{j=1} \lambda_j / n^2_j } } \right )^{\frac{2p}{p+2}} ,
}
which leads to
\bes{
\bbE \nm{F - \widehat{F}}_{L^p_{\mu}(\cX ; \cY)} \lesssim_p  \sqrt{\frac{n}{m}} \sigma + B^{\frac{p}{p+2}} L^{\frac{2}{p+2}} \left (d \sum^{d}_{j=1} \frac{\lambda_j}{n^2_j} \right )^{\frac{1}{p+2}} + L \sqrt{\sum_{j > d} \lambda_j}.
}
Now recall that $n = n_1 \cdots n_d$. Set $n_i = \lfloor c \sqrt{\lambda_i} \rfloor$ for some $c$ sufficiently large so that $c \sqrt{\lambda_i} \geq 1$, $\forall i = 1,\ldots,d$. Then, using the fact that $\lfloor x \rfloor \geq x/2$ for $x \geq 1$, we have
\bes{
\bbE \nm{F - \widehat{F}}_{L^p_{\mu}(\cX ; \cY)} \lesssim_p  \sqrt{\frac{c^d \sqrt{\lambda_1 \cdots \lambda_d}}{m }} \sigma + B^{\frac{p}{p+2}} L^{\frac{2}{p+2}} d^{\frac{2}{p+2}} c^{-\frac{2}{p+2}}  + L \sqrt{\sum_{j > d} \lambda_j}.
}
We now set
\bes{
c = \left ( \frac{m}{\sigma^2 \sqrt{\lambda_1 \cdots \lambda_d } } \right )^{\frac{p+2}{(p+2)d+4}} \left ( B^{2p} L^4 d^4 \right )^{\frac{1}{(p+2)d+4}} 
}
to obtain
\bes{
\bbE \nm{F - \widehat{F}}_{L^p_{\mu}(\cX ; \cY)} \lesssim_p  \left ( \frac{m}{\sigma^2 \sqrt{\lambda_1\cdots\lambda_d} } \right )^{-\frac{2}{(p+2)d+4} } (B^{2p} L^4 d^4 )^{\frac{d}{2((p+2)d+4)}} +L \sqrt{\sum_{j>d} \lambda_j} .
}
This holds, provided
\bes{
c \geq 1/\sqrt{\lambda_d},
}
which is equivalent to
\bes{
 \frac{m}{\sigma^2}   \geq \frac{(\lambda_1\cdots\lambda_d)^{\frac{1}{2}}}{\left ( B^{2p} L^4 d^4 \right )^{\frac{1}{p+2}} (\lambda_d)^{\frac{(p+2) d+4}{2(p+2)}} }, 
}
as required.
}

\prf{[Proof of Theorem \ref{t:main-ub-white}]
Arguing as in the previous proof, we deduce that
\bes{
\bbE \nm{F - \widehat{F}}_{L^p_{\mu}(\cX ; \cY)} \lesssim_p  B w^{-t}_{r} + \sqrt{\frac{r n}{m}} \sigma + B^{\frac{p}{p+2}} L^{\frac{2}{p+2}} \left (d \sum^{d}_{j=1} \frac{\lambda_j}{n^2_j} \right )^{\frac{1}{p+2}} + L \sqrt{\sum_{j > d} \lambda_j}.
}
We now proceed in an identical manner.
}

\section{Proofs of the results in Section \ref{ss:main-result}}\label{s:proofs-of-main-res}

\subsection{Lower bounds}\label{s:proofs-lbs}

\prf{
[Proof of Proposition \ref{prop:no-alg-conv}]
This is an immediate consequence of Theorem \ref{t:main-lb}.
}

\prf{[Proof of Theorems \ref{thm:main-exp-decay} and \ref{thm:main-exp-decay-white}; lower bounds]
For completeness, we establish a lower bound for $\betanew > 0$, as discussed in Remark \ref{rem:beta0-1}.
For convenience, let $\mu_d = \sum^{d}_{j=1} 1/\lambda_j$ and $s = L^2 m / (c \sigma^2)$. Using Theorem \ref{t:main-lb}, we can write
\be{
\label{Gsd-proof}
\begin{split}
\cL_{m}(\cF ; L^p_{\mu}) & \lesssim \log \left ( \frac{\sqrt{\mu_d}}{L}\right ) + \frac{d}{p} \log \left ( \frac{p+1}{\iota} \right ) + \max \left \{ \log \left( \frac{8}{a} \right ) , \frac{1}{2+d} \log \left ( \frac{s}{\mu_d} \right ) - \frac{d}{2+d} \log(a) \right \} 
\\
& = : G(s,d)
\end{split}
}
for all sufficiently large $d$. 
Now suppose that $s , d \rightarrow \infty$ with
\be{
\label{kd-cond}
\frac{\log(s) - \log(\mu_d)}{d} \sim \frac{\log(s)}{d} \text{ and } \frac{\log(s)}{d} \rightarrow \infty.
}
Then
\be{
\label{Gkd-twiddle}
G(s,d) \sim \frac12 \log(\mu_d)+ \frac{d}{p} \log \left ( \frac{p+1}{\iota} \right ) +\frac{\log(s)}{d}  .
}
We now treat the two cases $0 < \betanew \leq 1$ and $\betanew > 1$ separately. Suppose first that $0 < \betanew \leq 1$ and let
\be{
\label{d-sqrtk-choice}
d = d(s) = \left \lfloor\sqrt{\log(s)} \right \rfloor.
}
Notice that \ef{kd-cond} holds for this choice of $d$, since $\mu_{d} \leq d \exp(\alphanew d^{\betanew})$ and therefore $\log(\mu_{d(s)}) = \cO(\sqrt{\log(s)}) = o(\log(s))$ as $s \rightarrow \infty$. With this, we have 
\bes{
G(s,d(s)) \asymp_{\iota,p} \sqrt{\log(s)},\quad s \rightarrow \infty.
}
Writing $\log(s) = \log(m/\sigma^2) + \log(L^2 / c)$, we deduce that
\bes{
\cL_m(\cF_{B,L} ; L^p_{\mu}) \lesssim_{\iota,p} \sqrt{\log(m/\sigma^2)},\quad m \rightarrow \infty,
}
as required.

Now consider $\betanew > 1$. 
Observe that $\mu_d = \sum^{d}_{i=1} \exp(\alphanew i^\betanew)$ satisfies
\bes{ 
\exp(\alphanew d^{\betanew}) \leq \mu_d \leq d \exp(\alphanew d^{\betanew}).
}
Hence, since $\betanew > 1$,
\bes{
G(s,d) \sim \frac{\alphanew}{2} d^{\betanew} + \frac{\log(s)}{d} ,\quad s , d \rightarrow \infty.
}
We now set
\bes{
d(s) = \left \lfloor \left ( \log(s) \right )^{\frac{1}{\betanew+1}} \right \rfloor
}
and observe that \ef{kd-cond} holds for this choice. We deduce that 
\bes{
G(s,d(s)) \asymp_{\alphanew} (\log(s))^{\frac{\betanew}{\betanew+1}},
}
and therefore
\bes{
\cL_m(\cF_{B,L} ; L^p_{\mu}) \lesssim_{\alphanew} (\log(m/\sigma^2))^{\frac{\betanew}{\betanew+1}},\quad m \rightarrow \infty,
}
as required.
}

\prf{[Proof of Theorems \ref{thm:main-alg-decay} and \ref{thm:main-alg-decay-white}; lower bounds]
The argument is similar to the previous case. Let $\mu_d = \sum^{d}_{j=1} 1/\lambda_j$, $s = L^2 m / (c \sigma^2)$ and use Theorem \ref{t:main-lb} to write 
\eas{
\cL_{m}(\cF ; L^p_{\mu}) &\lesssim \log \left ( \frac{\sqrt{\mu_d}}{L}\right ) + \frac{d}{p} \log \left ( \frac{p+1}{\iota} \right ) + \max \left \{ \log \left(\frac{8}{a} \right) , \frac{1}{2+d} \log \left ( \frac{s}{\mu_d} \right ) - \frac{d}{2+d} \log(a) \right \} 
\\
& = : G(s,d)
}
for all sufficiently large $d$. Now suppose once more that \ef{kd-cond} holds, so that $G$ satisfies \ef{Gkd-twiddle}. Observe that $\mu_d = \sum^{d}_{j=1} j^{\alphanew} \sim \frac{d^{\alphanew+1}}{\alphanew+1}$, $d \rightarrow \infty$, and therefore
\bes{
G(s,d) \sim  \frac{d}{p} \log \left ( \frac{p+1}{\iota} \right ) + \frac{\log(s)}{d} .
}
We now set \ef{d-sqrtk-choice} once more and argue in the same way to obtain
\bes{
\cL_m(\cF_{B,L} ; L^p_{\mu}) \lesssim_{\iota,p} \sqrt{\log(m/\sigma^2)},\quad m \rightarrow \infty,
}
as required.
}

\rem{
[Non-sharpness of the lower bounds]
\label{rem:why-lower-bounds-nonsharp}
As noted, we currently lack tight characterizations of the minimax risk for exponentially-decaying eigenvalues with $0 < \betanew < 1$ (see Remark \ref{rem:beta0-1}) or algebraically-decaying eigenvalues for any $\alphanew > 1$. 
This stems from the factor $(\iota / (p+1))^d$
in \ref{t:main-lb}. Since $\iota \leq 1$ (see Assumption \ref{ass:mu}), this factor is always exponentially-small in $d$. When the eigenvalues decay algebraically or exponentially with $0 < \betanew < 1$, the decay of this term dominates in the asymptotic behaviour of the function $G(s,d)$. This is seen in the previous proofs. For exponential decay with $\betanew \geq 1$ this term no longer dominates, with the consequence being a tight characterization. We believe the factor $(\iota / (p+1))^d$ 
is an artefact of the proof of Theorem \ref{t:main-lb}. Removing it is an objective of future work.
}

\subsection{Upper bounds}\label{s:proofs-ubs}

\prf{[Proof of Theorem \ref{thm:main-exp-decay}; upper bound]
We may assume without loss of generality that $B , L \geq 1$, since this only potentially increases the size of the set $\cF_{B,L}$.
Suppose that $\lambda_i = \exp(-\alphanew i^{\betanew})$. Then
\bes{
\lambda_1 \cdots \lambda_d \leq \exp(-c_{\alphanew,\betanew} d^{\betanew+1}),
}
for some constant $c_{\alphanew,\betanew} > 0$
and
\eas{
\sum_{j > d} \lambda_j &\leq \int^{\infty}_{d} \exp(-\alphanew t^\betanew) \D t.
\\
& = \frac{1}{\betanew \alphanew^{1/\betanew}} \int^{\infty}_{\alphanew d^\betanew} u^{1/\betanew-1} \exp(-u) \D u
\\
& = \frac{1}{\betanew \alphanew^{1/\betanew}} \Gamma(1/\betanew , \alphanew d^\betanew),
}
where $\Gamma$ denotes the upper incomplete Gamma function. It is known that $\Gamma(a,z) \sim z^{a-1} \E^{-z}$ as $z \rightarrow \infty$. We deduce that
\bes{
\sum_{j > d} \lambda_j \lesssim_{\alphanew,\betanew} d^{1-\betanew} \exp(-\alphanew d^{\betanew}).
}
Let $k = m/ \sigma^2$. Then Theorem \ref{t:main-ub-hilbert} gives that
\bes{
\cM_{m}(\cF_{B,L} ; L^p_{\mu}) \lesssim_{\alphanew,\betanew,B,L,p} d^{\frac{4}{2(p+2)}} k^{-\frac{2}{(p+2)d+4}} \exp(-c_{\alphanew,\betanew,p} d^{\betanew}) + d^{1/2-\betanew/2} \exp(-c_{\alphanew,\betanew,p} d^{\betanew} )
}
for some constant $c_{\alphanew,\betanew,p} > 0$, provided
\bes{
 k \geq \exp(c'_{\alphanew,\betanew,p} d^{\betanew+1}) 
}
for some constant $c'_{\alphanew,\betanew,p} > 0$. Here we used the facts that $B,L \geq 1$ and $d \geq 1$.
Set
\bes{
d = \left \lfloor \left ( \frac{\log(k)}{c'_{\alphanew,\betanew,p}} \right )^{\frac{1}{\betanew+1}} \right \rfloor.
}
For sufficiently large $m$, we have $k,d \geq 1$ and therefore
\eas{
\cM_m(F_{B,L} ; L^p_{\mu}) & \lesssim d \exp(-c_{\alphanew,\betanew,p} d^\betanew) \lesssim_{\alphanew,\betanew,p} (\log(m/\sigma^2))^{\frac{1}{\betanew+1}}\exp \left ( -c''_{\alphanew,\betanew,p} (\log(m / \sigma^2))^{\frac{\betanew}{\betanew+1}} \right ),\quad m \rightarrow \infty .
}
Consequently,
\bes{
\cL_{m}(\cF_{B,L} ; L^p_{\mu}) \gtrsim_{\alphanew,\betanew,p} (\log(m / \sigma^2))^{\frac{\betanew}{\betanew+1}}, \quad m \to \infty,
}
as required.
}

\prf{
[Proof of Theorem \ref{thm:main-exp-decay-white}; upper bound]
We proceed as in the previous proof, using Theorem \ref{t:main-ub-white} in place of Theorem \ref{t:main-ub-hilbert} and with $k = m / (r \sigma^2)$. 
Choosing the same value for $d$ and noting that $w_r = r$ by assumption, we deduce that
\bes{
\cM_{m}(\cF_{B,L} ; L^p_{\mu}) \lesssim_{\alphanew,\betanew,B,L,p} r^{-t} + \exp \left ( -c''_{\alphanew,\betanew,p} (\log(k))^{\frac{\betanew}{\betanew+1}} \right ),\quad k \rightarrow \infty,
}
and therefore
\bes{
\cL_{m}(\cF^t_{B,L} ; L^p_{\mu}) \gtrsim_{\alphanew,\betanew,p} \min \left \{ (\log(k))^{\frac{\betanew}{\betanew+1}} , t \log(r)  \right \} =  \min \left \{ \left ( \log(m/\sigma^2) - \log(r) \right )^{\frac{\betanew}{\betanew+1}} , t \log(r) \right \}.
}
This holds for all $r \in \bbN$. Now let $r = r(m)$ the smallest integer such that
\bes{
 \log(r) \geq (\log(m/\sigma^2))^{\frac{\betanew}{\betanew+1}}.
}
Notice that $\log(r) = o(\log(m/\sigma^2))$ as $m \rightarrow \infty$. The result now follows.
}

\prf{
[Proof of Theorem \ref{thm:main-alg-decay}; upper bound]
We proceed in a similar manner. Let $\lambda_i = i^{-\alphanew}$ and notice that
\bes{
\lambda_1 \cdots \lambda_d = (d!)^{-\alphanew} \leq (d/\E)^{-\alphanew d}
}
and 
\bes{
\sum_{j > d} \lambda_j \lesssim_{\alphanew} d^{1-\alphanew}
}
Hence Theorem \ref{t:main-ub-hilbert} gives that
\bes{
\cM_{m}(\cF_{B,L} ; L^p_{\mu}) \lesssim_{\alphanew,B,L,p} d k^{-\frac{2}{(p+2)d+4}} \left ( \frac{d}{\E} \right )^{-\frac{\alphanew d}{(p+2)d+4} } + d^{\frac{1-\alphanew}{2}},
}
provided
\bes{
k \geq (d/\E)^{-\alphanew d/2} d^{\alphanew \frac{(p+2) d+4}{2(p+2)} - \frac{4}{p+2}}.
}
This is equivalent to
\be{
\label{k-d-cond-ub-alg}
\log(k) \geq \frac{\alphanew}{2} d + \frac{2(\alphanew-2)}{p+2} \log(d).
}
Set
\bes{
d = \left \lfloor \frac{4}{(p+2)+\alphanew p } \frac{\log(k)}{\log(\log(k))} \right \rfloor
}
and observe that \ef{k-d-cond-ub-alg} holds for all sufficiently large $k$. Also note that $d \to \infty$ as $k \to \infty$. Consequently, substituting this value into the above error bound, we deduce that
\bes{
\cM_{m}(\cF_{B,L} ; L^p_{\mu}) \lesssim_{\alphanew,B,L,p}
d k^{-\frac{2}{(p+2)d} } d^{-\frac{\alphanew}{p+2}} + d^{\frac{1-\alphanew}{2}}
\lesssim_{\alphanew, p} \left ( \frac{\log(k)}{\log(\log(k))} \right )^{\frac{1-\alphanew}{2}}, \quad k \to \infty,
}
as required.
}

\prf{
[Proof of Theorem \ref{thm:main-alg-decay-white}; upper bound]
We proceed as in the previous proof, to deduce that
\bes{
\cL_{m}(\cF^t_{B,L} ; L^p_{\mu}) \lesssim_{\alphanew,B,L,p} r^{-t} + \left ( \frac{\log(k)}{\log(\log(k))} \right )^{\frac{1-\alphanew}{2}},
}
as $k \rightarrow \infty$, where $k = m/(r \sigma^2)$. Therefore,
\bes{
\cM_{m}(\cF^t_{B,L} ; L^p_{\mu}) \geq \min \left \{ \frac{\alphanew-1}{2} \log \left ( \frac{\log(k)}{\log(\log(k))} \right ) , t \log(r) \right \} +  \cO_{\alphanew,B,L,p}(1).
}
Now let $r = r(m)$ be the smallest integer such that
\bes{
t \log(r) \leq \frac{\alphanew-1}{2} \log \left ( \frac{\log(m/\sigma^2)}{\log(\log(m/\sigma^2))} \right ).
}
Observe that $k \rightarrow \infty$ as $m \rightarrow \infty$ with this choice and also that
\bes{
\frac{\log(k)}{\log(\log(k))} \sim \frac{\log(m/\sigma^2)}{\log(\log(m/\sigma^2))},\quad m \rightarrow \infty.
}
We deduce that 
\bes{
\cL_m(\cF^t_{B,L} ; \cL^p_{\mu}) \geq \frac{\alphanew-1}{2} \log(\log(m/\sigma^2)) + \cO_{\alphanew,\betanew,B,L,p}(1),
}
as required.
}

\rem{[Other choices of $w_r$]
\label{wr-choice}
The above proofs of the upper bounds in Theorems \ref{thm:main-exp-decay-white} and \ref{thm:main-alg-decay-white} indicate how the choice of the weights $w_i$ in the smoothness spaces \ef{Yt-def} can influence the bounds. In general, we have that the log-minimax risk $\cL_m(\cF^t_{B,L} ; L^p_{\mu})$ in the white noise case is bounded below by (up to potential constants)
\bes{
\min \left \{ \log(w_r) , \cR(\log(m/\sigma^2) - \log(r)) \right \},
}
where $\cR : (0,\infty) \rightarrow (0,\infty)$ is some increasing function which provides the
analogous term in the lower bound for $\cL_m(\cF_{B,L} ; L^p_{\mu})$ in the case of Hilbert-valued Gaussian noise. Proposition \ref{prop:no-alg-conv} asserts that $\cR$ must grow sublinearly, i.e., $\cR(z) = o(z)$ as $z \rightarrow \infty$.
We are interested in when the lower bound in the white noise case matches that in the Hilbert-valued Gaussian noise case. This holds whenever it is possible to choose $r = r(m)$ such that
\be{
\label{r-cond-general}
\log(r) \leq c \log(m/\sigma^2) ,\quad \log(w_r) \gtrsim \cR(\log(m/\sigma^2)),
}
for some $0 < c < 1$ and all sufficiently large $m$.
This clearly holds when $w_i = i$, as considered in this work. Slower growth of $w_i$ is also possible, depending on behaviour of $\cR$. For example, if $\cR(z) \asymp \log(z)/\log(\log(z))$, 
as in the case of algebraically-decaying eigenvalues (see Theorem \ref{thm:main-alg-decay}), then the choice $w_r = \log(r)$ also satisfies \ef{r-cond-general}.
}

\section{Further discussion} \label{sec:further}

Before concluding, we now discuss several further topics.

\subsection{Finite-dimensional rates} \label{sec:finite-dimensional-rates}

Although it is not the primary focus of this paper, we now briefly show that our main upper and lower bounds yield the well-known minimax risks for finite-dimensional Lipschitz functions on bounded domains. To this end, let $\cX = \bbR^d$, $\cY = \bbR$, both equipped with the Euclidean inner product, and suppose that $\mu$ is the uniform probability measure on $[-1,1]^d$. Notice that $a = 1$, $b = 1/2$ and $\lambda_1 = \cdots = \lambda_d = 1$ in this case. 
Theorem \ref{t:main-lb} gives the lower bound
\bes{
\cM_{m}(\cF_{B,L} ; L^p_{\mu}) \gtrsim_{d,p,L,B} \left ( \frac{m}{\sigma^2} \right )^{-\frac{1}{2+d}},\qquad m \rightarrow \infty.
}
Technically speaking, this is only true if
\ef{d-cond-lb} holds, i.e., $B \geq L (8 \sqrt{d})$. However, this poses no issue thanks to Remark \ref{rem:ben-lazy}.

We now establish the matching upper bound. On the fact of it, Theorem \ref{t:main-ub-hilbert} would yield a suboptimal rate depending on $p$. This stems from the general types of measures $\mu$ allowed in this theorem. However, applying Lemma \ref{lem:fixed-design-histogram-bd} or \ref{lem:random-design-histogram-bd} with $R > 1$ and the choice $n_1 = \cdots = n_d = \lfloor n^{\frac1d} \rfloor$ immediately yields the upper bound
\bes{
\cM_m(\cF_{B,L} ; L^p_{\mu}) \lesssim_{d,p,L,B} \sqrt{\frac{n}{m}} \sigma + n^{-\frac1d},\qquad m \rightarrow \infty,
}
valid for all $n \in \bbN$. Optimizing with respect to $n$ we deduce that
\bes{
\cM_m(\cF_{B,L} ; L^p_{\mu}) \lesssim_{d,p,L,B} \left ( \frac{m}{\sigma^2} \right )^{-\frac{1}{2+d}},\qquad m \rightarrow \infty.
}
Combining with the lower bound above, we now recover the known minimax rate
\bes{
\cM_m(\cF_{B,L} ; L^p_{\mu}) \asymp_{d,p,L,B} \left ( \frac{m}{\sigma^2} \right )^{-\frac{1}{2+d}},\qquad m \rightarrow \infty.
}

\subsection{Double-exponential decay}\label{sec:double-exp-decay}

Suppose that $\lambda_i = \exp(- \exp(\alphanew i))$ has double-exponential decay. Then, similar to Theorems \ref{thm:main-exp-decay} and \ref{thm:main-exp-decay-white}, one can establish the following tight characterizations of the log-minimax risk.

\thm{[Tight characterization for double-exponentially decaying $\lambda_i$; Hilbert-valued Gaussian noise]
\label{thm:main-double-exp-decay}
Suppose that Assumption \ref{ass:mu} and Assumption \ref{ass:hilbert-noise} hold. Let $1 \leq p < \infty$ and suppose that $\lambda_i = \exp(-\exp(\alphanew i) )$ for some $\alphanew > 0$. Then
\bes{
\cL_m(\cF_{B,L} ;  L^p_{\mu})  \asymp_{\alphanew,\iota,p} \frac{\log(m / \sigma^2)}{\log(\log(m / \sigma^2))},\quad m \rightarrow \infty,
}
where $\cL_m(\cF_{B,L} ; L^p_{\mu}) = - \log(\cM_m(\cF_{B,L} ;  L^p_{\mu}))$ is as in \ef{neg-log-minimax}, with $\cM_m(\cF_{B,L} ;  L^p_{\mu})$ given by either  \ef{minimax-fixed} or \ef{minimax-random}.
}

\thm{[Tight characterization for double-exponentially decaying $\lambda_i$; white noise]
\label{thm:main-double-exp-decay-white}
Suppose that Assumption \ref{ass:mu} and Assumption \ref{ass:white-noise} hold. Let $t > 0$, $1 \leq p < \infty$ and suppose that $w_i = i$ and $\lambda_i = \exp(-\exp(\alphanew i) )$ for some $\alphanew > 0$. Then
\bes{
\cL_m(\cF^t_{B,L} ; , L^p_{\mu})  \asymp_{\alphanew,\iota,t,p} \frac{\log(m / \sigma^2)}{\log(\log(m / \sigma^2))},\quad m \rightarrow \infty.
}
where $\cL_m(\cF_{B,L} ; L^p_{\mu})$ is as in \ef{neg-log-minimax}, with $\cM_m(\cF_{B,L} ; , L^p_{\mu})$ given by either \ef{minimax-fixed} or \ef{minimax-random}.
}

These results mean that algebraic decay rates, while impossible due to Proposition \ref{prop:no-alg-conv}, can `nearly' be achieved for sufficiently fast-decaying eigenvalues. Indeed, they imply that
\bes{
\cM_{m}(\cF ; L^p_{\mu}) \asymp C_1 \left ( \frac{m}{\sigma^2} \right )^{-\frac{C_2}{\log(\log(m/\sigma^2))} },\quad m \rightarrow \infty,
}
for constants $C_1,C_2 > 0$ depending on the various parameters. Thus, the rate is algebraic for a double exponentially-large range of $m/\sigma^2$: specifically, for any $c > 0$ the error decays algebraically at least as fast as $(m/\sigma^2)^{-C_2/c}$ for all $m$ satisfying $m/\sigma^2 \leq \exp(\exp(c))$.

\prf{[Proofs of Theorems \ref{thm:main-double-exp-decay} and \ref{thm:main-double-exp-decay-white}]
For the lower bounds for the minimax risk, we argue as in Section \ref{s:proofs-lbs}. Let $\mu_d = \sum^{d}_{j=1} 1/\lambda_j$ and $s = L^2 m / (c \sigma^2)$ once more and $G(s,d)$ be as in \ef{Gsd-proof}.
Observe that $\log(\mu_d) \sim \exp(\tau d)$ as $d \rightarrow \infty$. Thus, assuming that \ef{kd-cond} holds, we have
\bes{
G(s,d) \sim \frac12 \exp(\tau d) + \frac{\log(s)}{d},
}
as $s , d \rightarrow \infty$. Now set
\bes{
d = d(s) = \left \lfloor \frac{1}{\tau} \log \left ( \frac{ \log(s)}{\log( \log(s))} \right ) \right \rfloor
}
and observe that \ef{kd-cond} holds in this case. This gives
\bes{
G(s,d) \asymp_{a} \frac{\log(s)}{\log(\log(s))}.
}
The lower bound now follows after writing $\log(s) = \log(m/\sigma^2) + \log(L^2/c)$.

For the upper bound, we proceed as in Section \ref{s:proofs-ubs}. We first observe that
\bes{
\lambda_1 \cdots \lambda_d \leq \exp(-\exp(a d)).
}
Also,
\bes{
\sum_{i > d} \lambda_i = \sum^{\infty}_{j=1} \exp(-\exp(a d) \exp(a j)) \leq \sum^{\infty}_{j = 1} \exp(-a \exp(a d) j ) \leq 2 \exp(-a \exp(a d)),
}
for all sufficiently large $d$. Let $k = m/\sigma^2$. Then  Theorem \ref{t:main-ub-hilbert} implies that
\be{
\label{Mm-upper-bd-double-exp}
\cM_{m}(\cF_{B,L} ; L^p_{\mu})  \lesssim_{p,B,L} k^{-\frac{2}{(p+2)d+4}} \exp \left ( - \frac{\exp(a d)}{(p+2) d + 4} \right ) d^{\frac{2}{p+2}} + \exp \left ( - \frac{a \exp(a d)}{2} \right ),
}
provided
\be{
\label{k-upper-cond-double-exp}
k \geq  (B^{2p} L^4)^{- \frac{1}{2 p}}  \exp \left ( \frac{(p+2)(d-1)+4}{2(p+2)}\exp(a d) \right ).
}
Let 
\bes{
d = \left \lfloor \frac{1}{a} \log \left ( \frac{a \log(k)}{ \log(\log(k))} \right )  \right \rfloor 
}
and observe that
\bes{
\frac{(p+2)(d-1)+4}{2(p+2)} \exp(a d) \sim \frac{1}{2 a} \log(\log(k)) \frac{a \log(k)}{\log(\log(k))} = \frac12 \log(k).
}
Therefore \ef{k-upper-cond-double-exp} holds for all large $k$. We now return to \ef{Mm-upper-bd-double-exp}. With this choice of $d$, observe that
\bes{
k^{-\frac{2}{(p+2)d + 4}} = \exp \left ( - \frac{2\log(k)}{(p+2) d+4} \right ) \leq \exp \left(-c_{p,a} \frac{\log(k)}{\log(\log(k))} \right ).
}
Also,
\bes{
\exp \left ( - \frac{\exp(a d)}{(p+2) d + 4} \right ) d^{\frac{2}{p+2}} \lesssim \exp \left ( - c_{p,a} \frac{\log(k)}{(\log(\log(k)))^2} \right ) \left( \frac{\log(\log(k)}{a} \right)^{\frac{2}{p+2}} \lesssim 1,
}
for a potentially different value of $c_{p,a}$.
Using this, we deduce that
\bes{
\cM_{m}(\cF_{B,L} ; L^p_{\mu}) \lesssim_{p,B,L} \exp \left(-c_{p,a} \frac{\log(k)}{\log(\log(k))} \right ),\quad k \rightarrow \infty.
}
This gives the result in the case of Hilbert-valued Gaussian noise. For Gaussian white noise we argue similarly, using Theorem \ref{t:main-ub-white} instead of Theorem \ref{t:main-ub-hilbert}.
}

\subsection{Minimax risk bounds for $C^{k,\alpha}$ operators}

In this final part, we consider smoother classes of operators. Our goal is to demonstrate that higher regularity does not improve the minimax rate, up to potential constants.

Let $k \in \bbN_0$ and $0 < \alpha \leq 1$. Let $\cL^k(\cX ; \cY)$ denote the space of continuous, $k$-linear maps $\cX^{\otimes k} \rightarrow \cY$. Given an operator $F : \cX \rightarrow \cY$, we write $D^k F(X) \in \cL^k(\cX ; \cY)$ for its Fr\'echet derivative at $X \in \cX$, whenever it exists. We then define $C^{k,\alpha}(\cX ; \cY)$ as the space of $k$-times differentiable operators $\cX \rightarrow \cY$ whose $k$th derivative is H\"older-$\alpha$ continuous,
and equip it with the norm
\bes{
\nm{F}_{C^{k,\alpha}} = \sum^{k}_{j=0} \sup_{X \in \cX} \| D^j F(X) \|_{\cL^j(\cX ; \cY)} + [D^k F]_{\alpha} 
}
where
\bes{
[D^k F]_{\alpha} = \sup_{\substack{X , X' \in \cX \\ X \neq X'}} \frac{ \nm{D^k F(X) - D^k F(X')}_{\cL^k(\cX ; \cY)} }{\nm{X - X'}_{\cX}}.
}
For succinctness, we consider only Hilbert-valued Gaussian noise in this section (one could also consider Gaussian white noise using the same ideas as those employed earlier in this paper).
Given $B \geq 0$, we define the class
\be{
\label{F-k-alpha-def}
\cF^{k,\alpha}_{B} = \left \{ F \in C^{k,\alpha}(\cX ; \cY) : \nm{F}_{C^{k,\alpha}} \leq B \right \}.
}
Notice that $F^{0,1}_{B}$ is the class of Lipschitz operators with $\nm{F}_{C^{0,1}} \leq B $, which is essentially the same as that studied in the rest of this paper (the main difference is the use of only one constant $B$ instead of both $B$ and $L$).

We now present the following result, which is a generalization of the lower bound for Lipschitz operators (Theorem \ref{t:main-lb}) to $C^{k,\alpha}$ operators.
\thm{
\label{t:lb-Ckalpha}
Suppose that Assumptions \ref{ass:mu} and \ref{ass:hilbert-noise} hold and that $\lambda_1 \leq 1$. Let $k \in \bbN_0$, $0 < \alpha \leq 1$ and $\cF = \cF^{k,\alpha}_{B}$ be given by \ef{F-k-alpha-def}. Then there are constants $c = c_{\Upsilon} > 0$, $c_{k,\alpha} > 0$ and $c_p > 0$ such that
\be{
\label{main-fixed-lb-Ckalpha}
\cM_m(\cF ; L^p_{\mu}) \gtrsim \frac{B}{(c_{k,\alpha})^d \Lambda_{k,\alpha,d} } (\iota c_p )^{d/p} \left ( \min \left \{ 1, \frac{a}{8} , \left ( \frac{B^2 m}{c a^d \sigma^2 (c_{k,\alpha})^{2d} \Lambda^2_{k,\alpha,d}} \right )^{-\frac{1}{2(k+\alpha)+d}}  \right \} \right )^{k+\alpha}
}
for all $d \in \bbN$, 
where
\be{
\label{Lambda-kalphad-def}
\Lambda_{k,\alpha,d} = (k + 2) \left ( \sum^{d}_{i=1} (\lambda_i)^{-\frac{\alpha}{2-\alpha}} \right )^{1-\frac{\alpha}{2}} \left ( \sum^{d}_{j=1} 1/\lambda_j \right )^{\frac{k}{2}} 
}
and $\cM_{m} (\cF ; L^p_{\mu})$ is given by either \ef{minimax-fixed} or \ef{minimax-random}.
}
Note that the assumption $\lambda_1 \leq 1$ is made for convenience, as it simplifies the expression for $\Lambda_{k,\alpha,d}$. Note also that there is no condition on $d$, unlike in Theorem \ref{t:main-lb}. This is because of the use of only one constant $B$ in the definition of $\cF$, as per the discussion above. 
The proof follows the same lines as that of Theorem \ref{t:main-lb}, and commences with the construction of a well-separated set (Lemma \ref{l:well-separated}).

\lem{
[Construction of a well-separated set of $C^{k,\alpha}$ operators]
\label{l:well-separated-Ckalpha}
Let $k \in \bbN_0$, $0 < \alpha \leq 1$, $0 < h \leq \min \{1, a/8 \}$ and $1 \leq p \leq \infty$, where $a$ is as in Assumption \ref{ass:mu}. Suppose that $\lambda_1 \leq 1$. Then there are constants $c_{k,\alpha} > 0$ and $c_p > 0$ and a collection of $C^{k,\alpha}$ functionals $F_0,\ldots,F_M : \cX \rightarrow \bbR$ such that
\be{
\label{Fj-Ckalpha-bounds}
\begin{split}
\nm{F_j}_{C^{k,\alpha}} & \leq B,\quad 0 \leq j \leq M,
\\
|F_j(X)| & \leq \frac{B h^{k+\alpha}}{(c_{k,\alpha})^d \Lambda_{k,\alpha,d} },\quad \forall X \in \cX,\ 0 \leq j \leq M,
\end{split}
}
and
\be{
\label{Fjk-lower-kalpha}
\nm{F_j - F_{k}}_{L^p_{\mu}(\cX ; \bbR)} \gtrsim \frac{B }{(c_{k,\alpha})^d \Lambda_{k,\alpha,d}} (\iota c_{p})^{d/p}  h^{k+\alpha},\quad 0 \leq j < k \leq M,
}
where $M$ satisfies $\log(M) \geq (a/h)^{d} \log(2) / 8$ and
$\Lambda_{k,\alpha,d}$ is as in \ef{Lambda-kalphad-def}.
In the case $p = \infty$, the exponent $d/p$ in \ef{Fjk-lower-kalpha} is interpreted as $0$.
}
\prf{
We once more break the proof into several steps.

\pbk
\textit{Step 1a: Single-bump construction and norm estimation.} 
Let $u : \bbR \rightarrow [0,\infty)$ be a $C^{\infty}$-bump function, with $1 = u(0) \geq u(x)$, $\forall x$, and $u(x) = 0$ for $|x| \geq 1$.  Let $f : \bbR^d \rightarrow \bbR$ be defined by
\bes{
f(x) = \prod^{d}_{i=1} u(x_i / \sqrt{\lambda_i}),\quad x \in \bbR^d,
}
and $F : \cX \rightarrow \bbR$ be defined by
\bes{
F(X) = f(\ip{X}{\phi_1}_{\cX},\ldots,\ip{X}{\phi_d}_{\cX}),
}
where the $\lambda_i$ and $\phi_i$ are as in Assumption \ref{ass:mu}. We first need to estimate the $C^{k,\alpha}$-norm of $F$. For this, observe that, for any $l = 0,\ldots,k$ and $X,H_1,\ldots,H_l \in \cX$,
\be{
\label{Dk-F-equiv-f}
D^l F(X)(H_1,\ldots,H_l) = D^l f(x)(h_1,\ldots,h_l),
}
where $x = (\ip{X}{\phi_j}_{\cX})^{d}_{j=1}$, $h_i = (\ip{H_i}{\phi_j}_{\cX})^{d}_{j=1}$. This follows from the chain rule or Fr\'echet derivatives, after noticing that $F$ is the composition of $f$ with the linear operator $\cX \rightarrow \bbR^d, X \mapsto ( \ip{X}{\phi_j}_{\cX} )^{d}_{j=1}$. This implies that
\bes{
\nm{D^l F(X)}_{\cL^l(\cX ; \bbR)} = \nm{D^l f(x)}_{\cL^l(\bbR^d ; \bbR)} \leq \sqrt{\sum^{d}_{j_1,\ldots,j_k = 1} | \partial_{j_1} \cdots \partial_{j_l} f(x) |^2 }.
}
Now observe that
\bes{
\partial_{j_1} \cdots \partial_{j_l} f(x) = (\lambda_{j_1} \cdots \lambda_{j_l})^{-1/2} \prod^{d}_{i=1} u^{(r_i)}(x_i / \sqrt{\lambda_i}) ,
}
where $r_i \in \bbN_0$ with $r_i \leq l$. Now let $c_k = \nm{u}_{C^k}$ so that
\bes{
|\partial_{j_1} \cdots \partial_{j_l} f(x) | \leq (\lambda_{j_1} \cdots \lambda_{j_l})^{-1/2} (c_k)^d.
}
We deduce that
\be{
\label{Dk-F-bound}
\nm{D^l F(X)}_{L^k(\cX ; \bbR)} \leq (c_k)^d \left ( \sum^{d}_{j=1} 1/\lambda_j \right )^{\frac{l}{2}} \leq (c_k)^d \left ( \sum^{d}_{j=1} 1/\lambda_j \right )^{\frac{k}{2}},\quad l =0,\ldots,k.
}
Here, in the inequality, we used the fact that $\lambda_1 \leq 1$.
We also need to estimate the semi-norm $[D^k F]_{\alpha}$. It follows from \ef{Dk-F-equiv-f} and Parseval's identity that
\bes{
[D^k F]_{\alpha} = [D^k f]_{\alpha}.
}
To estimate the latter, we let $x,x' \in \bbR^d$ and write
\be{
\label{Dkf-partial-split}
\begin{split}
&\nm{D^k f(x) - D^k f(x')}_{\cL^k(\bbR^d ; \bbR)} 
\\
&= \sup_{\substack{h_1,\ldots,h_k \in \bbR^d \\ \nm{h_1}_2 = \cdots = \nm{h_k}_2 = 1}} \left | \sum^{d}_{j_1,\ldots,j_k = 1} \left ( \partial_{j_1} \cdots \partial_{j_k} f(x) - \partial_{j_1} \cdots \partial_{j_k} f(x')  \right ) (h_1)_{j_1} \cdots (h_k)_{j_k} \right |
\\
& \leq \left [ \sum^{d}_{j_1,\ldots,j_k = 1}  | \partial_{j_1} \cdots \partial_{j_k} f(x) - \partial_{j_1} \cdots \partial_{j_k} f(x') |^2 \right ]^{\frac12} .
\end{split}
}
Now write
\be{
\label{f-lambda-g}
\partial_{j_1} \cdots \partial_{j_k} f(x) - \partial_{j_1} \cdots \partial_{j_k} f(x') = (\lambda_{j_1} \cdots \lambda_{j_k})^{-1/2} (g(x) - g(x')),
}
where
\bes{
g(x) = \prod^{d}_{i=1} u^{(r_i)}(x_i/\sqrt{\lambda_i})
}
for $r_i \in \bbN_0$ with $r_i \leq k$.
Using a telescoping sum, we write
\eas{
g(x) - g(x')  =&\  [g(x_1,\ldots,x_d) - g(x_1,\ldots,x_{d-1},x'_d) ]
\\
& + [g(x_1,\ldots,x_{d-1},x'_d) - g(x_1,\ldots,x_{d-2},x'_{d-1},x'_d) ] + \cdots ,
}
to obtain
\bes{
|g(x) - g(x')| \leq \sum^{d}_{i=1} [ u^{(r_i)} ]_{\alpha} (c_k)^{d-1} |x_i - x'_i |^{\alpha} / (\lambda_i)^{\alpha/2}.
}
Now, recall that $u$ is supported in $[-1,1]$. Hence, if $r_i < k$, then we have
\bes{
[u^{(r_i)}]_{\alpha} \leq 2 \nm{u^{(r_i+1)}}_{C^0} \leq 2 \nm{u}_{C^{k,\alpha}}
}
and, trivially, when $r_i = k$,
\bes{
[u^{(r_i)}]_{\alpha} \leq \nm{u}_{C^{k,\alpha}}.
}
Therefore,
\bes{
|g(x) - g(x')| \leq (c_{k,\alpha})^{d} \sum^{d}_{i=1} |x_i - x'_i |^{\alpha} / (\lambda_i)^{\alpha/2},
}
where $c_{k,\alpha} : = 2\nm{u}_{C^{k,\alpha}} \geq c_k$. Using H\"older's inequality, 
\bes{
|g(x) - g(x')| \leq (c_{k,\alpha})^{d} \nm{x - x'}^{\alpha}_2 \left ( \sum^{d}_{i=1} (\lambda_i)^{-\frac{\alpha}{2-\alpha}} \right )^{1-\frac{\alpha}{2}}.
}
Combining this with \ef{Dkf-partial-split} and \ef{f-lambda-g} we deduce that
\eas{
\nm{D^k f(x) - D^k f(x')}_{\cL^k(\bbR^d ; \bbR)} 
 \leq (c_{k,\alpha})^d \nm{x - x'}^{\alpha}_2 \left ( \sum^{d}_{i=1} (\lambda_i)^{-\frac{\alpha}{2-\alpha}} \right )^{1-\frac{\alpha}{2}} \left ( \sum^{d}_{j=1} 1/\lambda_j \right )^{\frac{k}{2}}.
}
We deduce that
\bes{
[D^k F ]_{\alpha} = [D^k f ]_{\alpha} \leq (c_{k,\alpha})^d  \left ( \sum^{d}_{i=1} (\lambda_i)^{-\frac{\alpha}{2-\alpha}} \right )^{1-\frac{\alpha}{2}} \left ( \sum^{d}_{j=1} 1/\lambda_j \right )^{\frac{k}{2}}.
}
Combining this with \ef{Dk-F-bound} and using the fact that $\lambda_1 \leq 1$ by assumption, we conclude that
\bes{
\nm{F}_{C^{k,\alpha}} \leq (k + 2) (c_{k,\alpha})^d  \left ( \sum^{d}_{i=1} (\lambda_i)^{-\frac{\alpha}{2-\alpha}} \right )^{1-\frac{\alpha}{2}} \left ( \sum^{d}_{j=1} 1/\lambda_j \right )^{\frac{k}{2}} .
}
Since $(k+2) \leq (k+2)^d$, after replacing $c_{k,\alpha}$ by $(k+2) c_{k,\alpha}$ we now deduce that
\bes{
\nm{F}_{C^{k,\alpha}} \leq (c_{k,\alpha})^d \Lambda_{k,\alpha,d}.
}

\pbk
\textit{Step 1b: Bump function construction.} We now proceed as in the proof of Lemma \ref{l:well-separated}. Let $C_1,\ldots,C_n \in \cX$, $\theta \in \{0,1\}^n$ and define
\bes{
F_{\theta}(\cdot) = \frac{B h^{k+\alpha}}{(c_{k,\alpha})^d \Lambda_{k,\alpha,d}} \sum^{n}_{i=1} \theta_i F \left ( \frac{\cdot - C_i}{h} \right ),
}
so that \ef{Fj-Ckalpha-bounds} holds. 

\pbk
\textit{Step 2: Estimating the $L^2_{\mu}$-distance of two bump functions.} The argument is once more similar. For $1 \leq p < \infty$, we have
\bes{
\nm{F_{\theta} - F_{\theta'} }^p_{L^p_{\mu}(\cX ; \bbR)} = \left ( \frac{B h^{k+\alpha}}{(c_{k,\alpha})^d \Lambda_{k,\alpha,d} } \right )^p \sum^{n}_{i=1} |\theta_i - \theta'_i | I_{i,p},
}
where
\bes{
I_{i,p} = \int_{\cX} \left | F \left ( \frac{X - C_i}{h} \right ) \right |^p \D \mu(X).
}
We now write
\bes{
I_{i,p} \geq (b h)^d \left ( \int^{1}_{-1} u(x)^p \D x \right )^d = : h^d (2 c_{p} b )^d.
}
where $c_{p} = \nm{u}^p_{L^p(\bbR)}/2$.  We deduce that
\bes{
\nm{F_{\theta}-F_{\theta'}}_{L^p_{\mu}(\cX ; \bbR)} \geq \frac{B }{(c_{k,\alpha})^d \Lambda_{k,\alpha,d}} (2 c_{p} b)^{d/p}  h^{k+\alpha+d/p} H(\theta,\theta')^{1/p}
}
and similarly for $p = \infty$.

\pbk
\textit{Steps 3 and 4.} These are identical to the corresponding steps in the proof of Lemma \ref{l:well-separated}. Hence we omit the details.
}

\prf{[Proof of Theorem \ref{t:lb-Ckalpha}]
The proof follows the same lines as that of Theorem \ref{t:main-lb}. Let $Y = Y_1$ be the eigenvector of the covariance operator of $\Upsilon$ corresponding to its largest eigenvalue $\upsilon_1$ and define $G_j(X) = Y F_j(X)$ for all $X \in \cX$ and $j = 0,\ldots,M$, where the $F_j$ are as in Lemma \ref{l:well-separated-Ckalpha}. By construction, we have $G_j \in \cF^{k,\alpha}_{B}$. We now let $d(\cdot,\cdot)$ denote the $L^p_{\mu}$-distance, $\kappa = 1/16$ and
\bes{
s^* = \frac{c_0}{2} \frac{B}{(c_{k,\alpha})^d \Lambda_{k,\alpha,d}} (\iota c_p)^{d/p} h^{k+\alpha},
}
where $c_0>0$ is the universal constant from \ef{Fjk-lower-kalpha} and the terms $c_{k,\alpha}$, $c_p$ and $\Lambda_{k,\alpha,d}$ are as in Lemma \ref{l:well-separated-Ckalpha}. Once more, we have $\log(M) \geq \log(2)$. In the fixed design case, we have
\bes{
D(P_j \| P_0) = \frac{1}{2 \sigma^2} \sum^{m}_{i=1} \nm{\Upsilon^{-1/2} G_j(X_i) }^2_{\cY} \leq \frac{m}{2 \sigma^2 \upsilon_1} \max_{X \in \cX} | F_j(X) |^2.
}
We now apply \ef{Fj-Ckalpha-bounds} to get that
\bes{
D(P_j \| P_0) \leq c \frac{m B^2 h^{2(k+\alpha)}}{\sigma^2 (c_{k,\alpha})^{2d} \Lambda^2_{k,\alpha,d}},
}
where $c = c_{\Upsilon} > 0$. The same bound also applies in the random design case. We now apply Theorem \ref{t:fano} along with the fact that $\log(M) \geq (a/h)^d \log(2)/8$ to deduce that
\bes{
\cM_m(\cF ; L^p_{\mu}) \gtrsim \frac{B}{(c_{k,\alpha})^d \Lambda_{k,\alpha,d}} (\iota c_p)^{d/p} h^{k+\alpha},
}
provided
\bes{
\frac{m B^2 h^{2(k+\alpha)}}{\sigma^2 (c_{k,\alpha})^{2d} \Lambda^2_{k,\alpha,d}} \leq c (a/h)^d,
}
for a potentially different constant $c = c_{\Upsilon}$. We now set
\bes{
h = \min \left \{ 1 , \frac{a}{8} , \left ( \frac{B^2 m}{c a^d \sigma^2 (c_{k,\alpha})^{2d} \Lambda^2_{k,\alpha,d}} \right )^{-\frac{1}{2(k+\alpha)+d}}  \right \}
}
and substitute this into the previous expression to get the result.
}

Using this result, we immediately deduce the following, which is a generalization of Proposition \ref{prop:no-alg-conv}. It confirms that algebraic decay is impossible for \textit{any} finite H\"older regularity.

\prop{[Algebraic decay is impossible for $C^{k,\alpha}$ operators]
Suppose that Assumptions \ref{ass:mu} and \ref{ass:hilbert-noise} hold. Let $k \in \bbN_0$, $0 < \alpha \leq 1$ and $1 \leq p < \infty$. Then, for any $q > 0$, we have
\bes{
\limsup_{m \rightarrow \infty}  \cM_{m}(\cF ; L^p_{\mu}) \cdot m^q  = + \infty,
}
where $\cF = \cF^{k,\alpha}_{B}$.
}

We next present the following generalization of Theorem \ref{thm:main-exp-decay}, which gives a tight characterization of the minimax rate in the case of exponentially-decaying $\lambda_i$.

\thm{[Tight characterization for exponentially-decaying $\lambda_i$; $C^{k,\alpha}$ operators with Hilbert-valued Gaussian noise]
Suppose that Assumption \ref{ass:mu} and Assumption \ref{ass:hilbert-noise} hold. Let $k \in \bbN$, $0 < \alpha \leq 1$ and $1 \leq p < \infty$ and suppose that $\lambda_i = \exp(-\alphanew i^{\betanew})$ for some $\alphanew > 0$ and $\betanew \geq 1$. Then
\bes{
\cL_m(\cF^{k,\alpha}_{B} ;  L^p_{\mu})  \asymp_{k,\alpha,\alphanew,\betanew,\iota,p} \left(\log(m / \sigma^2) \right)^{\frac{\betanew}{\betanew+1}},\quad m \rightarrow \infty,
}
where $\cL_m(\cF^{k,\alpha}_{B} ; L^p_{\mu}) = - \log(\cM_m(\cF^{k,\alpha}_{B} ;  L^p_{\mu}))$ is as in \ef{neg-log-minimax}, with $\cM_m(\cF^{k,\alpha}_{B} ;  L^p_{\mu})$ given by either \ef{minimax-fixed} or \ef{minimax-random}.
}

Upon comparison with the Lipschitz case (Theorem \ref{thm:main-exp-decay}), the main conclusion of this result is that the minimax rate does not change with increasing smoothness, except possibly for the constants. Note that we have made no attempt to track the dependence on $k$.

\prf{
To lower bound $\cL_m(\cF^{k,\alpha}_{B} ;  L^p_{\mu})$, we simply use the lower bound from Theorem \ref{thm:main-exp-decay}, since $\cF^{k,\alpha}_{B} \subseteq \cF_{B,B}$ for $k \in \bbN$. We now consider the upper bound and proceed as in Section \ref{s:proofs-lbs}. Once more, we consider arbitrary $\omega > 0$ in the proof.
Let $s = B^2 m / (c \sigma^2)$. Then, using Theorem \ref{t:lb-Ckalpha}, we have
\eas{
\cL_m(\cF ; L^p_{\mu}) \lesssim & \ \log \left ( \frac{\Lambda_{k,\alpha,d}}{B} \right ) + d \left ( \log(c_{k,\alpha}) - \frac{\log(\iota c_p)}{p} \right ) 
\\
& + (k+\alpha) \max \left \{ 0, \log \left ( \frac{8}{\alphanew} \right ) , \frac{1}{2(k+\alpha)+d} \log \left ( \frac{s}{(c_{k,\alpha})^{2d} \Lambda^2_{k,\alpha,d}} \right ) - \frac{d}{2(k+\alpha)+d} \log(a) \right \}
\\
= : & \ G(s,d)
}
Now suppose that $s,d \rightarrow \infty$ with
\be{
\label{asymp-ass-Ckalpha}
\frac{1}{2(k+\alpha)+d} \log \left ( \frac{s}{(c_{k,\alpha})^{2d} \Lambda^2_{k,\alpha,d}} \right ) \sim \frac{\log(s)}{d}\text{ and } \frac{\log(s)}{d} \rightarrow \infty.
}
Then
\bes{
G(s,d) \sim \log(\Lambda_{k,\alpha,d}) + d \left ( \log(c_{k,\alpha}) - \frac{\log(\iota c_p)}{p} \right ) + (k+\alpha) \frac{\log(s)}{d}.
}
Suppose first that $0 < \betanew \leq 1$ and set
\bes{
d = d(s) = \left \lfloor \sqrt{\log(s)} \right \rfloor.
}
Notice that
\bes{
\Lambda_{k,\alpha,d} \leq  d^{1+\frac{k-\alpha}{2}} \exp \left ( c_{k,\alpha,\alphanew}  d^{\betanew} \right )
}
and therefore
\bes{
\log(\Lambda_{k,\alpha,d}) = o (\log(s)),\quad s \rightarrow \infty.
}
We deduce that \ef{asymp-ass-Ckalpha} holds. Hence
\bes{
G(s,d) \asymp_{k,\alpha,\iota,p} \sqrt{\log(s)},\quad s \rightarrow \infty.
}
The result now follows for $0 < \betanew \leq 1$ after writing $\log(s) = \log(m/\sigma^2) + \log(B^2/c)$.

For $\betanew > 1$, we note that
\bes{
\exp \left ( c_{k,\alpha,\alphanew}  d^{\betanew} \right ) \leq \Lambda_{k,\alpha,d} \leq  d^{1+\frac{k-\alpha}{2}} \exp \left ( c_{k,\alpha,\alphanew}  d^{\betanew} \right )
}
and therefore
\bes{
G(s,d) \sim c_{k,\alpha,\alphanew} d^{\betanew} + (k+\alpha) \frac{\log(s)}{d}.
}
We now set $d(s) = \left \lfloor (\log(s))^{\frac{1}{\betanew+1}} \right \rfloor$ and observe that \ef{asymp-ass-Ckalpha} holds for this choice. Therefore
\bes{
G(s,d) \asymp_{k,\alpha,\alphanew} (\log(s))^{\frac{\betanew}{\betanew+1}},
}
as required.
}

For succinctness, we do not consider the case of algebraically- or double exponentially-decaying eigenvalues. However, we remark in passing that analogous versions of Theorems \ref{thm:main-alg-decay} and \ref{thm:main-double-exp-decay} can be established in these cases, up to possible changes in the constants depending on $k$, $\alpha$.

\section{Conclusion} \label{sec:conclusion}

We developed a minimax theory for operator learning from noisy input-output data when both the inputs and outputs may be infinite dimensional. Our general upper and lower bounds cover fixed and random designs and accommodate both Hilbert-valued Gaussian noise and Gaussian white noise. Our results make precise how statistical difficulty is governed by the eigenvalues of the covariance of $\mu$. Our results also recover classical finite-dimensional rates as a special case.

There are several open problems that remain. In the case of Lipschitz operators, most importantly, our lower bounds are not sharp for algebraically decaying eigenvalues and, more generally, for exponential decay with $0 < \betanew < 1$. A key technical obstacle is the $d$-dependence in the current lower-bound argument. Removing it would likely lead to a sharp characterization in these regimes. In the algebraic case, we conjecture that the true minimax rate is polylogarithmic in the sample size, with an exponent determined by the decay rate of the eigenvalues. Even in regimes where our bounds match, such as exponential and double-exponential eigenvalue decay, identifying the optimal constant in the exponent appears nontrivial. Finally, for very rapidly decaying eigenvalues (such as double-exponential decay) one can achieve `nearly' algebraic performance. Beyond Lipschitz operators, we showed that higher-order (H\"older) regularity does not change the bounds. This shows that imposing higher, but still finite, regularity cannot overcome the curse of sample complexity in operator learning. Extending this framework to other operator classes (e.g., holomorphic, where lower bounds are not known, or Besov classes, where neither upper nor lower bounds are known) is another natural direction for future work.

\section*{Acknowledgements}

BA acknowledges support from the Natural Sciences and Engineering Research Council of Canada
(NSERC) through grant RGPIN/2470-2021. 
GM acknowledges support from the Hausdorff Center for Mathematics (HCM)
in Bonn, funded by the Deutsche Forschungsgemeinschaft (DFG, German Research
Foundation) under Germany’s Excellence Strategy – EXC-2047/1 – 390685813.

\bibliographystyle{abbrv}
\bibliography{minimax-refs}

\end{document}